\documentclass{article}
\usepackage[a4paper]{geometry}
\usepackage[utf8]{inputenc}         
\usepackage{graphicx,longtable,booktabs,picture}
\usepackage{multirow}
\usepackage{colortbl}
\usepackage[dvipsnames,table,xcdraw]{xcolor}
\usepackage{hhline}
\usepackage{import}
\usepackage{makecell}
\usepackage[squaren,Gray]{SIunits}

\providecommand{\keywords}[1]
{
  \small	
  \textbf{\textit{Keywords---}} #1
}

\usepackage{tikz}
\tikzset{every path/.append style={line width=1pt}}
\usepackage{pgfplots}
\usepgfplotslibrary{fillbetween}
\usepgfplotslibrary{groupplots}
\pgfplotsset{compat=1.7,
	emphasize/.code args={#1:#2with#3}{
		\pgfplotsextra{
			\draw[fill=#3,fill opacity=0.5] ({axis cs:#1,0} |- {axis description cs:0,0}) 
			rectangle ({axis cs:#2,0} |- {axis description cs:0,1});
		}
	}
}
\pgfplotsset{
	bar group size/.style 2 args={
		/pgf/bar shift={%
			-0.5*(#2*\pgfplotbarwidth + (#2-1)*\pgfkeysvalueof{/pgfplots/bar group skip})  + 
			(.5+#1)*\pgfplotbarwidth + #1*\pgfkeysvalueof{/pgfplots/bar group skip}},%
	},bar group skip/.initial=2pt}

\usepackage[margin=20pt,font=small,labelfont={bf,sf},format=hang]{caption}
\usepackage{listings,multicol}
\usepackage{tabularx}
\usepackage{amsmath} 
 
\usepackage{amssymb}
\usepackage{cancel}
\usepackage{wasysym}
\usepackage{units} 
\usepackage{pgfplots}

\newtheorem{problem}{Problem}
\newtheorem{remark}{Remark}
\definecolor{maroon}{rgb}{0.69, 0.19, 0.38}
\usepackage{pdfpages}
\usepackage{lineno}
\modulolinenumbers[5]
\usepackage{hyperref}

\begin{document}

\title{A mixed phase-field fracture model\\ for crack propagation in punctured EPDM strips}

\author{Katrin Mang$^{1,}$\footnote{Corresponding author: \url{mang@ifam.uni-hannover.de}}, Andreas Fehse$^{2}$, Nils Hendrik Kr\"{o}ger$^{2,3}$, Thomas Wick$^{1}$\\
\small $^{1}$Institute of Applied Mathematics, Leibniz Universit{\"a}t Hannover, Welfengarten 1, 30167 Hannover, Germany\\
\small $^{2}$Deutsches Institut für Kautschuktechnologie e. V., Eupener Straße 33, 30519 Hannover, Germany\\
\small $^{3}$material prediction GmbH, Nordkamp 14, 26203 Wardenburg, Germany\\
}

\maketitle

%

\begin{abstract}
In this work, we present crack propagation experiments evaluated by
digital image correlation (DIC) for a carbon black filled ethylene propylene diene monomer rubber (EPDM) and numerical modeling with the help of variational phase-field fracture.
Our main focus is the evolution of cracks in one-sided notched EPDM strips containing a circular hole.
The crack propagation experiments are complemented with investigations identifying the mechanical material properties as well as the critical strain energy release rate.
For simulating the evolution of cracks with a given notch, phase-field fracture modeling is a popular approach. To avoid volume-locking effects considering fractures in nearly incompressible materials, a quasi-static phase-field
fracture model in its classical formulation is reformulated with the help of a mixed form of the solid-displacement equation.
The new established mixed phase-field fracture model is applied to simulate crack propagation in punctured EPDM strips by using the experimentally identified material parameters with mixed finite elements.
To discuss agreements and point out challenges and differences, the crack paths, the maximal force response, the traverse displacement at the crack start, as well as force-displacement curves of the experimental and numerical results are compared.
\end{abstract}



\section{Introduction}

In the recent past, many authors contributed to improve the understanding of the fatigue behavior of elastomeric materials. Applied fatigue models often concentrate on the end-of-life behaviour, e.g.\,by W\"{o}hler's method,  e.g.\,\cite{Wohler.1870,Gehrmann.2019}. Later approaches combine the classic Paris-Erdogan approach for crack propagation with statistical particle distribution and particle size distribution, e.g.\,\cite{Ludwig.2017,ElYaagoubi2018,Gehrmann.2019a}. Another alternative is the continuum damage mechanics approach, e.g.\,\cite{Lemaitre.1985,Grandcoin.2014}. 
A challenging task within modeling fatigue is the correct prediction of the crack propagation via finite element simulations. Classical approaches to simulate the growth of a crack are based on energetic failure criteria using the virtual crack extension method \cite{Charrier2003} or crack tip closure method \cite{Charrier2003,Timbrell2003}. Similar approaches consider cohesive elements, see e.\,g.\,\cite{Kaliske2014} for an application using inelastic cohesive models.\\
In mechanics, phase-field models for quasi-static brittle fracture based on Griffith's theory \cite{griffith1920phenomena} are a well-established approach to simulate complex crack phenomena as crack nucleation, propagation, branching or merging.
Phase-field modeling is based on a smoothed indicator variable called crack phase-field parameter describing a fracture path continuously over the displacement field. 
A first variational approach for quasi-static brittle fracture was introduced by Francfort and Marigo \cite{FraMar98}, see also \cite{bourdin2008variational}. A regularized formulation was presented two years later by Bourdin et al.~\cite{bourdin2000numerical}. This formulation is motivated by the regularization of Ambrosio and Tortorelli \cite{ambrosio1990approximation} for elliptic free-discontinuity problems. Often in the literature, 
the Ambrosio-Tortorelli functional with a quadratic energy degradation function (AT$_2$) is considered in the context of phase-field fracture problems. To decrease the impact of a possibly material-dependent length scale parameter $\epsilon$, Wu \cite{wu2017unified,wu2018geometrically,wu2018length} introduced a unified phase-field theory with the idea that the unknown phase-field and its gradient are introduced to regularize the sharp crack topology in a purely geometric context. See \cite{wu2017unified} for a detailed description. 
Furthermore, Amor et al. \cite{amor2009regularized} proposed a volumetric-deviatoric decomposition of the elastic energy density, because the regularized formulation does not distinguish between fracture behaviour in tension and compression \cite{ambati2015review}.
Detailed overviews on phase-field fracture modeling from mechanical and mathematical perspectives are given by Ambati et al. \cite{ambati2015review}, Wu et al. \cite{wu2018phase}, Bourdin and Francfort \cite{bourdin2019past}, and Wick \cite{wick2020multiphysics}.\\
If the considered solid is nearly incompressible, e.g. EPDM, the classical Euler-Lagrange equations derived from a regularized energy functional fail due to volume locking effects. Locking in finite element simulations leads to the problem, that the solid displacements are underestimated and unrealistic numerical results can be observed.
To allow also simulating crack growth in rubber-like materials, the phase-field fracture model is extended, see \cite{mang2020phase,mang2020pamm} for a detailed derivation. This approach builds on a mixed form of the solid-displacement equation resulting in two unknowns: 
a displacement field and a pressure variable. Via Taylor-Hood elements we achieve a stable discretization of the displacement-pressure system. Alternative approaches have been presented by Loew et al. \cite{Loew2019a,Loew2019b} for rate-dependent phase field damage models or by Faye et al. \cite{Faye2019} for mass sink models based on \cite{Volohk2017}.\\
In this work, the mixed phase-field fracture model is tested on punctured EPDM strips with a notch on 
one side, which are stretched until total failure. 
Being well aware of geometric and material non-linearities dealing with rubber and high strains, in this contribution we neglect those effects focusing on the modeling approach as well as on the qualitative comparison to the experiments. Whereas the experiments themselves are a benchmark for future studies. 
In addition, it is sufficient to assume plane-stress since the experimental specimens are thin and for this reason, two-dimensional numerical simulations are justified. With these simplifications 
we cannot expect perfect matches between numerical simulations and experiments. However, due to the complexity of the experimental setup and our newly developed
mixed phase-field fracture model using Wu's functional, such qualitative comparisons are 
a major effort.

To this end, the main contributions and the outline of this work are:
\begin{itemize}
 \item Determining the material compounding and properties as well as the (critical) strain energy release rate via digital image correlation (DIC) in Section \ref{Material} and \ref{Measure};
 \item Deriving a quasi-static mixed phase-field fracture model based on Wu's model \cite{wu2017unified} and Amor's strain energy splitting \cite{amor2009regularized} for incompressible solids similar to \cite{mang2020phase}, in Section \ref{notation} and \ref{mixedModel};
 \item Applying and substantiating the new model via numerical simulations of crack propagation in punctured EPDM strips in Section \ref{Simulation} and discussions 
 in Section \ref{Discussion}.
\end{itemize}
Section \ref{Conclusion} summarizes the content and the main results of this work.


\section{Experiments - mechanical characterization and crack path analysis}\label{Experi}


\subsection{Material compounding and sample preparation}\label{Material}
The experimental study is conducted using a sulphur crosslinked EPDM (Keltan $2450$) filled with $60$ phr carbon black N$550$. The EPDM mixture, see Table \ref{Tab:Material}, is prepared on a $5.0$ liter mixer adding the ingredients step-wise. The crosslinking agent and catalysts are admixed at a roller at 60$^\circ$C. The samples, $2 \mathrm{mm}$ thick specimens for all experiments, are compression moulded at 170$^\circ$C for six minutes corresponding to ``$t_{90}$'' plus two minutes of an according vulcameter test\footnote{``$t_{90}$'' denotes the time where the torque of the vulcameter curve reaches 90\% of its maximum.}. If required for the experiments, notches and circular inclusions are pierced in the specimens after vulcanisation.

\begin{table}[htbp!]
	\centering
	\begin{tabular}{|>{\columncolor[HTML]{C0C0C0}}c|c|c|c|c|c|c|c|c|c|}
		\hline
		Ingredients & \begin{tabular}[c]{@{}c@{}}EPDM \end{tabular} & \begin{tabular}[c]{@{}c@{}}Carbon black \\N$550$\end{tabular} & \begin{tabular}[c]{@{}c@{}}PEG- \\4000\end{tabular}  & \begin{tabular}[c]{@{}c@{}}Oil \end{tabular}  & \begin{tabular}[c]{@{}c@{}}Zinc  \\oxid\end{tabular}  & \begin{tabular}[c]{@{}c@{}}Stearic \\acid\end{tabular}  & Sulphur  & TBBS  & TBzTD  \\ \hline
		Phr& 100 &60 & 5.0 & 5.0 & 5.0 & 3.0 & 0.7 & 1.0  & 3.5 \\ \hline
		 \begin{tabular}[c]{@{}c@{}}Admixed at \end{tabular} & \multicolumn{6}{c|}{Mixer} & \multicolumn{3}{c|}{Roller} \\ \hline
		\begin{tabular}[c]{@{}c@{}}Admixed after \\ min\end{tabular} & 0 &\begin{tabular}[c]{@{}c@{}}2/3@1 \\ 1/3@2\end{tabular}& 2m30 & 1 & 2m30 & 1 & 1 & 1  & 1 \\ \hline
	\end{tabular}
	\caption{Recipe of the investigated EPDM mixture (phr $\widehat{=}$ parts per hundred parts of rubber related to mass parts).}
	\label{Tab:Material}
\end{table}


\subsection{Determination of mechanical material properties}\label{Measure}
\begin{figure}[htbp!]
	\centering
		\includegraphics[]{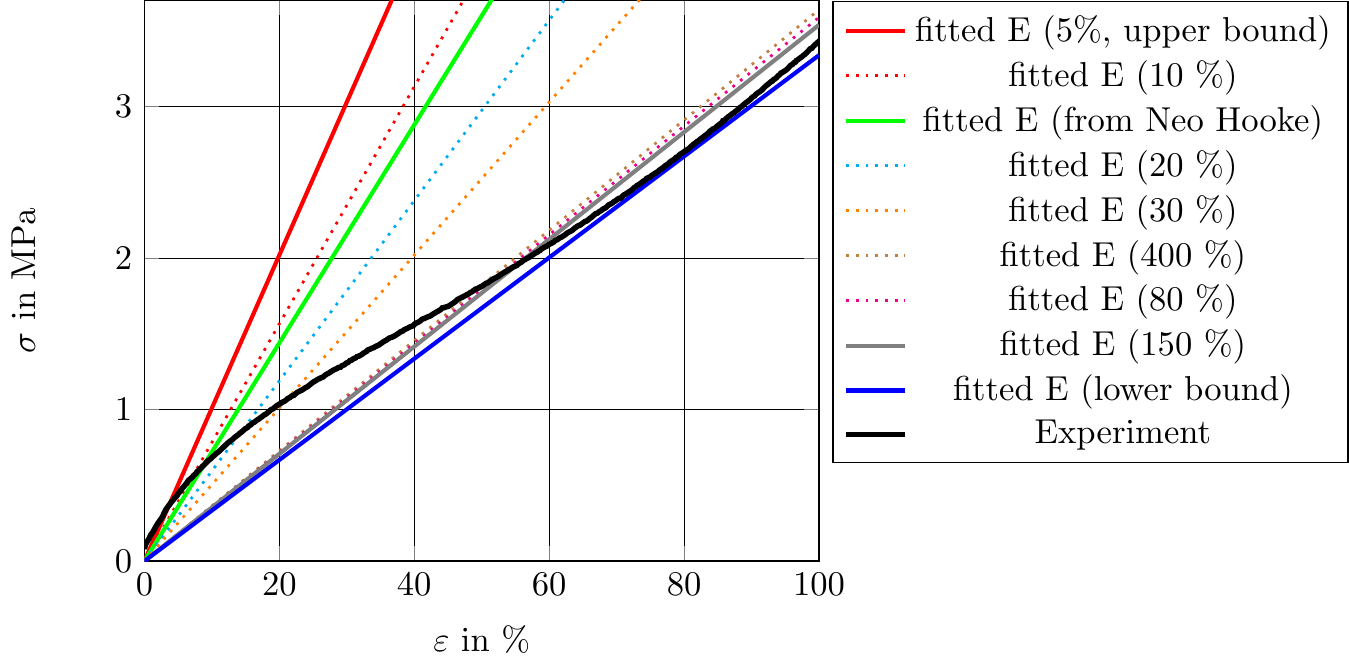}
		\caption{Uniaxial tension test (stress $\sigma$ versus strain $\varepsilon$) and prediction of the elastic behaviour based on parameter identifications of the elastic moduli related to specific fitting strain ranges.}
	\label{fig:EModul}
\end{figure}

\begin{figure}[htbp!]
	\centering
		\includegraphics[]{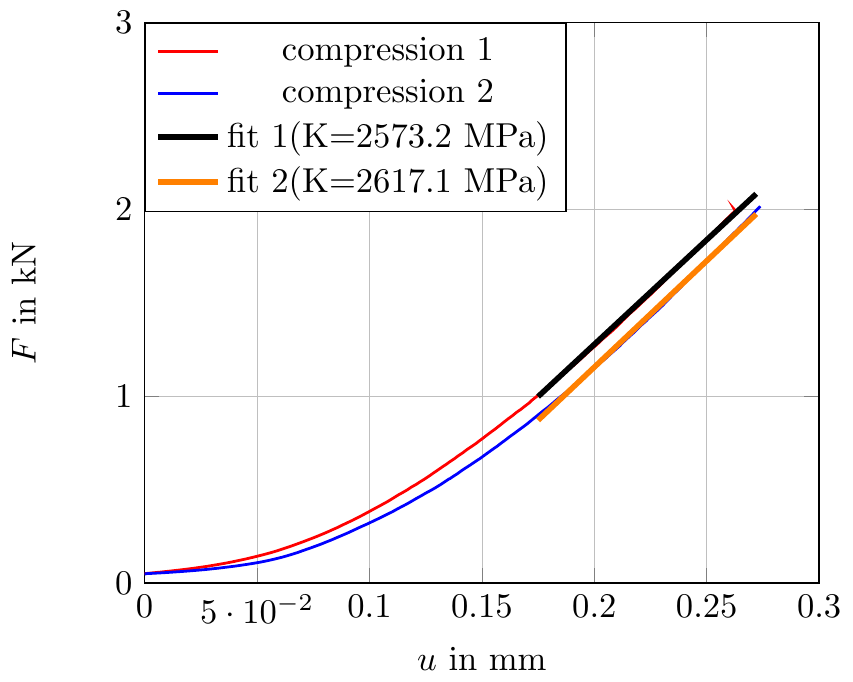}
	\caption{Volumetric compression test and estimation of the bulk modulus by a linear fit in the range between $1$ and $2 \mathrm{kN}$. The force $F$ and displacement $u$ are depicted in sense of compression, compared to tension the sign would switch. Please refer to \cite{Gehrmann2017,Ricker2019} for details on the experimental setup.}
	\label{fig:compression}
\end{figure}

\textbf{Elastic constants.} The mechanical material properties are determined on a ZWICK Universal test machine in four loading modes - uniaxial, planar (pure shear) and biaxial tension as well as confined (volumetric) compression. The tests were conducted with a (initial) strain rate of approximately $50\%$ per minute and a pre-force of $1 \mathrm{N}$, cf. Figure \ref{fig:EModul}(black line). The later defined material model uses a linear elastic approach such that the minor non-linear behaviour has to be estimated by e.g. the average behaviour until a certain strain, a lower or upper bound or alternative with an estimate for the Young's modulus in consistency with the Neo Hooke model (estimate holds for the initial stiffness, therefore one possibility for an upper bound definition), see Figure \ref{fig:EModul} (gray line for 'up to $150\%$ estimate' and green line for 'Neo Hooke estimate'). The bulk modulus is derived by volumetric compression experiments within the constant response region of $1000$ to $2000 \mathrm{N}$, see Figure \ref{fig:compression}. The parameter identification is done by using a minimization for the least square sum of the absolute errors between experiment and prediction considering the varying strain range. The resulting parameters are depicted in Tabular \ref{Tab:Youngsmodulus}. The estimation of linear material behavior might deviate for strains under 50\% but is a solid assumption for higher strains appearing in the later crack path experiments for this special EPDM. \\

\begin{table}[htbp!]
	\centering
\begin{tabular}{|>{\columncolor[HTML]{C0C0C0}}c|c|>{\columncolor[HTML]{ADD8E6}}c|c|c|c|>{\columncolor[HTML]{ADD8E6}}c|c|c|}
\hline
Parameter & upper  & Neo Hooke & 10 & 20 & 80 & 150 & 400 & lower  \\ \hline
$E$ [MPa] & 10.1  & 7.20 & 7.83 & 5.95  & 3.59 & 3.54 & 3.64 & 3.34 \\ \hline
$K$ [MPa] & 2595 & 2595 & 2595 & 2595 & 2595 & 2595 & 2595 & 2595\\ \hline 
$\lambda$ [MPa] & 2592.75  & 2593.40 & 2593.26 & 2593.68 & 2594.20 & 2594.21 & 2594.19 & 2594.26\\ \hline
$\mu$ [MPa] & 3.37 & 2.40 & 2.61& 1.98 & 1.20 & 1.18 & 1.21 & 1.11\\ \hline
$\nu$ [ ] & 0.49978 & 0.49985 & 0.49983 & 0.49987 & 0.49992 & 0.49992 & 0.49992 & 0.49993\\ \hline
\end{tabular}
	\caption{Young's modulus $E$, bulk modulus $K$, Lamé coefficients $\lambda$ and $\mu$ and Poisson ratio $\nu$ for different identification variants of $E$ (upper and lower bound, resp. results for parameter fit up to $X \%$ of strain, and estimated by Neo Hookean model), cf. Figures \ref{fig:EModul} and \ref{fig:compression}. Parameters set in light blue are used in the simulations in Section \ref{Simulation}.}
	\label{Tab:Youngsmodulus}
\end{table}

\begin{figure}[htbp!]
\centering
\includegraphics{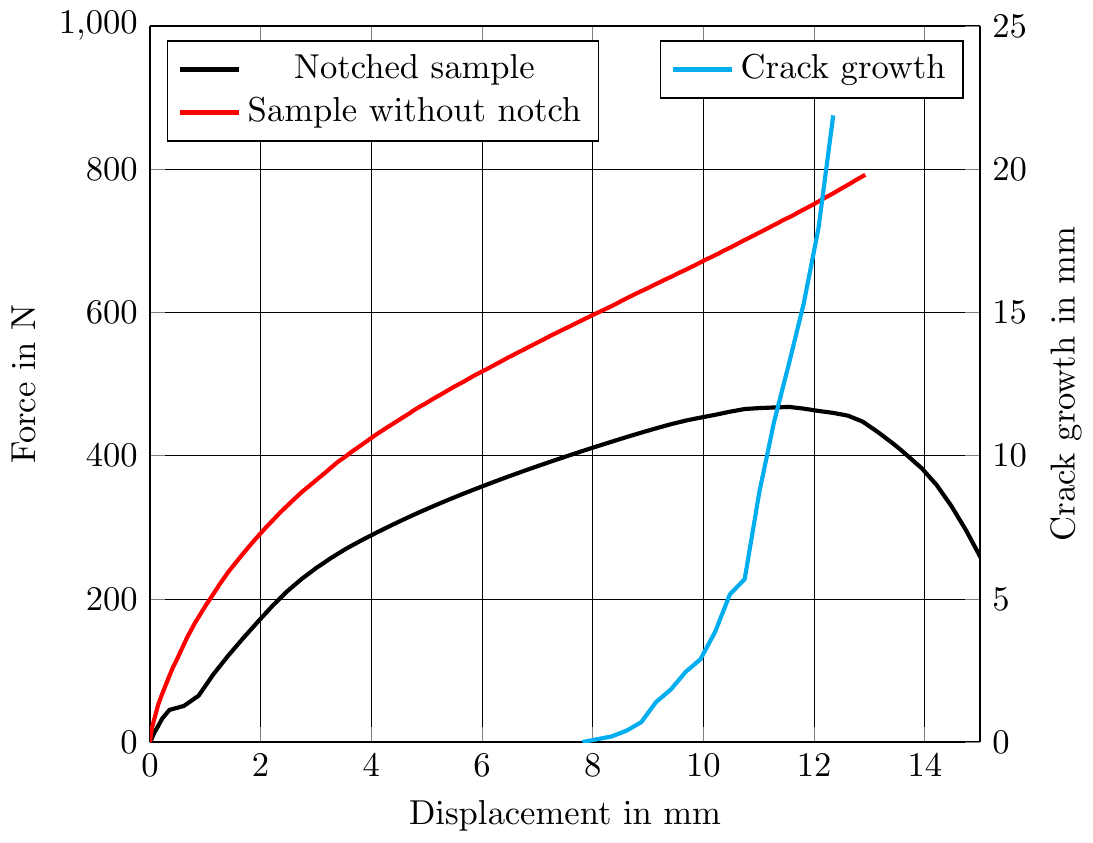}
\caption{Force-displacement curves (displacement $u$ in $y$-direction measured on the top boundary versus force $F_y$) for a notched sample and a sample without notch. Further, the crack growth evaluated from DIC, cf.\,Figure \ref{fig:tearingenergyspecimen}, in dependence of the displacement is plotted in blue.}
\label{fig:tearingenergyspecimen2}
\end{figure}

\textbf{Critical energy release rate.} 
Pure shear tests were used to determine the (critical) energy release rate $G_c$. The specimens had a length of $196 \mathrm{mm}$, a height of $28 \mathrm{mm}$ and a thickness of $1.8 \mathrm{mm}$. A ZWICK 1445 universal test machine was used for the tests. The traverse velocity was chosen to be $200 \mathrm{mm/min}$. Two types of tests were conducted: standard pure shear tests until break, and pure shear tests with notched samples. 
The notched sample had an initial crack $c_0$ of $47 \mathrm{mm}$.
In the second case, a square pattern was applied to the specimen such that the crack growth could be followed by a DIC recording. The crack growth was correlated with the recorded force-displacement curves. The resulting force-displacement and crack growth curves are shown in Figure \ref{fig:tearingenergyspecimen2}. The critical (strain) energy release rate can be determined using this data and the following formula from \cite{Roucou1,Roucou2}:

\begin{align*}
G_c=-\frac{\partial U}{\partial A}\Biggr|_\lambda = -\frac{d U}{d A}\Biggr|_\lambda.
\end{align*}

Here, $U$ is the stored elastic energy and $dU$ is the difference between the sample with and without a crack, \cite{Roucou1,Roucou2}. With $F_w$ and $u_w$ being the force-displacement data of the sample without a crack, $F_c$ and $u_c$ being the data of the sample with a crack, and $\lambda$ as the stretch, $dU$ reads as
\begin{align*}
dU=\int_\lambda F_w du_w-\int_\lambda F_c du_c.
\end{align*}

\begin{figure}[htbp!]
\centering
\includegraphics[width=0.6\textwidth]{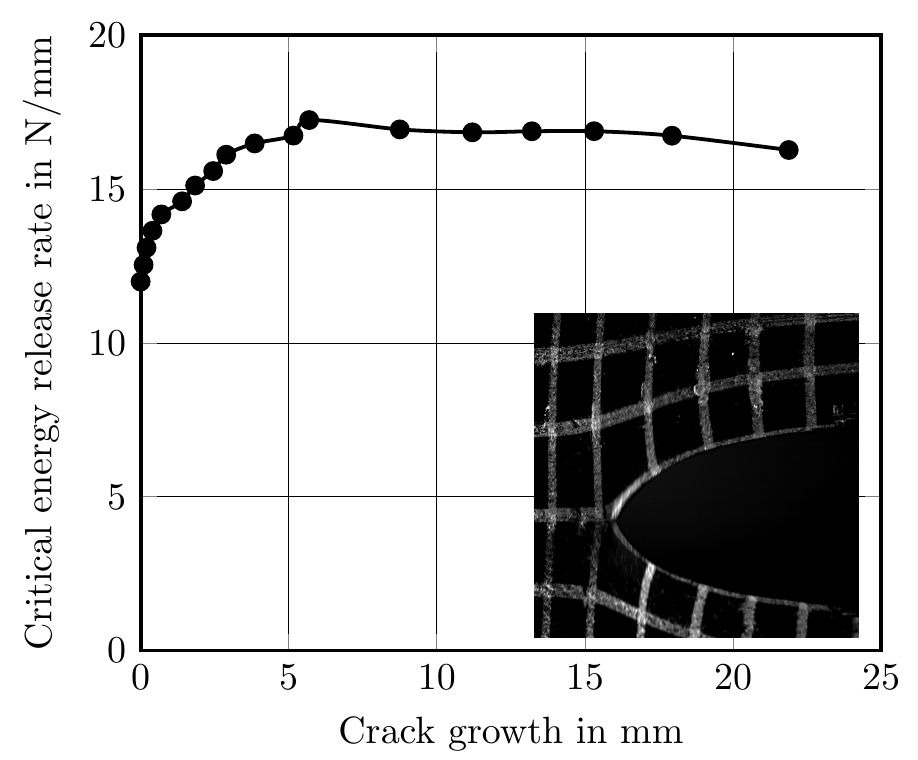}

\caption{Estimated (critical) strain energy release rate $G_c$ in dependence of the crack length of a single-edge notched pure shear tension test. The average value in the stable crack growth region after $5.0\mathrm{mm}$ crack length is approximately $17.0 \mathrm{N/mm}$.}
\label{fig:tearingenergyspecimen}
\end{figure}

The crack area $dA$ depends on the specimen thickness $t$ and the crack length $c_0 + c$. Here, $c_0$ is the initial crack length of $47 \mathrm{mm}$ and $c$ is the length of the growing crack, cf.\,\cite{Roucou1,Roucou2}:

\begin{align*}
dA= (c_0+c) \cdot t.
\end{align*}
Evaluating the experimental results the stationary value of the critical strain energy release rate $G_c$ can be approximated by $17.0\mathrm{N/mm}$, see Figure \ref{fig:tearingenergyspecimen}. The result coincides in general with literature results, see e.\,g.\,\cite{Roucou2019E} for carbon black filled styrene butadiene rubber.

\subsection{Crack path experiments}\label{Setup}
\begin{figure}[htbp!]
	\centering
		\reflectbox{\includegraphics[width=0.9\textwidth]{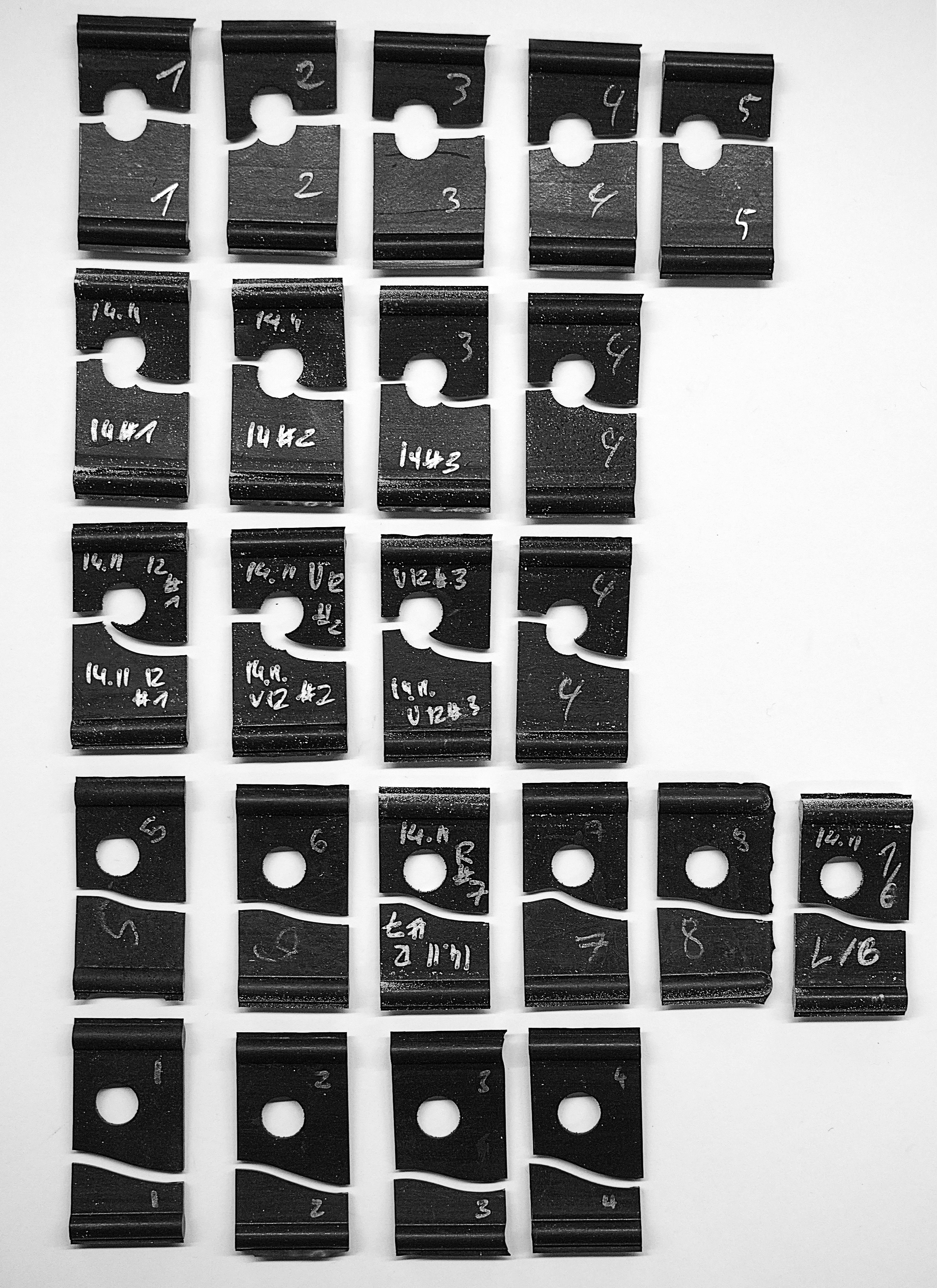}}
	\caption{Overview of the executed tests concerning tracking of the crack paths, back side of tested punctured strips, with an initially given notch of $1\mathrm{mm}$ length and $18$, $14$, $12$, $10$ and $6 \mathrm{mm}$ notch height (top to bottom).}
	\label{fig:NHK00169_v2}
\end{figure}

\begin{figure}[htbp!]
	\centering
		\includegraphics[width=0.85\textwidth]{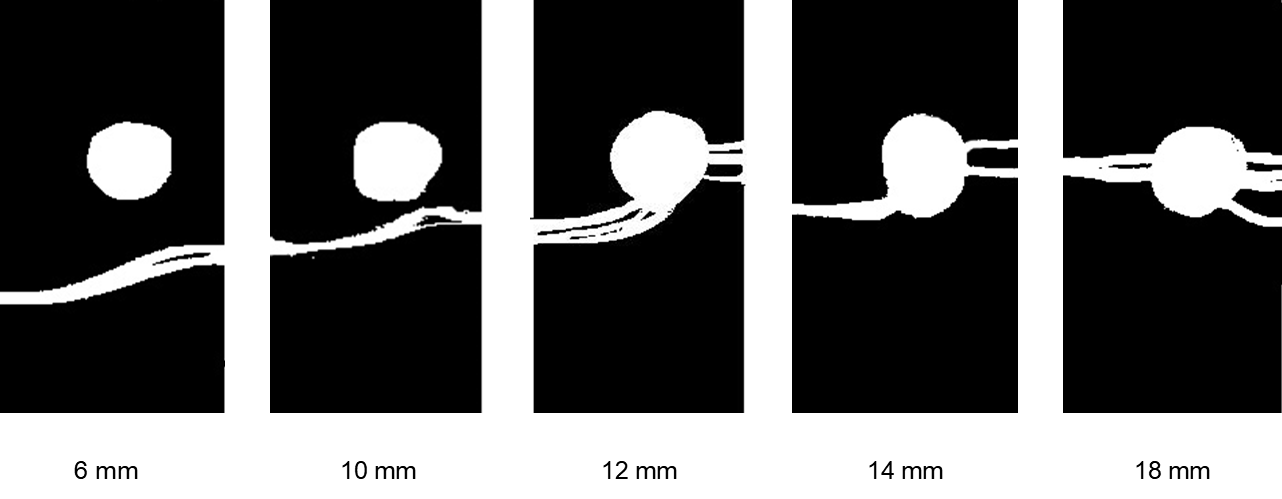}
	\caption{Evaluation of the averaged crack paths ($4$ to $6$ experiments, cf. Figure \ref{fig:NHK00169_v2}) of the tested punctured EPDM strips with given notches at a height of $6, 10, 12, 14$ and $18 \mathrm{mm}$ measured from the bottom boundary above the bottom bulges (left to right).}
	\label{fig:crackpaths}
\end{figure}

\begin{figure}[htbp!]
\centering
\begin{tikzpicture}
\begin{axis}[
 width=11cm,
    height=6cm,
    ylabel = $u$ at force maximum in $\mathrm{mm}$,
    xlabel = Notch height in $\mathrm{mm}$,
    grid =major,
    ybar,
    symbolic x coords={6,10,12,14,18},xtick={6,10,12,14,18},
    ymin=0,
    legend pos=outer north east,
]
\addplot[fill=red,opacity=0.75] coordinates
{(6,21.26454) (10,21.06456) (12,20.79792) (14,22.9977) (18,23.8876)};
\addlegendentry{Sim Neo Hooke, geo ref};
\addplot[fill=gray,opacity=0.75,error bars/.cd, y dir=both, y explicit] coordinates
{(6,25.0249) -=(0,1.0314) +=(0,2.5317) (10,27.2577) -=(0,4.4085) +=(0,4.1686) (12,25.6027) -=(0,2.5972) +=(0,4.7123) (14,26.3190) -=(0,4.3881) +=(0,7.1868) (18,22.6466) -=(0,3.8702) +=(0,4.4707)};
\addlegendentry{$4$ to $6$ Experiments};
\addplot[fill=blue,opacity=0.75] coordinates
{(6,29.6637) (10,27.29727) (12,27.26394) (14,32.46342) (18,34.06326)};
\addlegendentry{Sim $150\%$ fit, geo ref};
\end{axis}
\end{tikzpicture}
\caption{Traverse displacement $u$ in $y$-direction at the first force maximum of experiments and finite element simulations. In gray (middle bar) the average displacement of ($4$ to $6$) experiments,
in red and blue (left and right bar for each notch height), the numerically achieved displacements at the first force maximum are given.}
\label{fig:crackstart}
\end{figure}
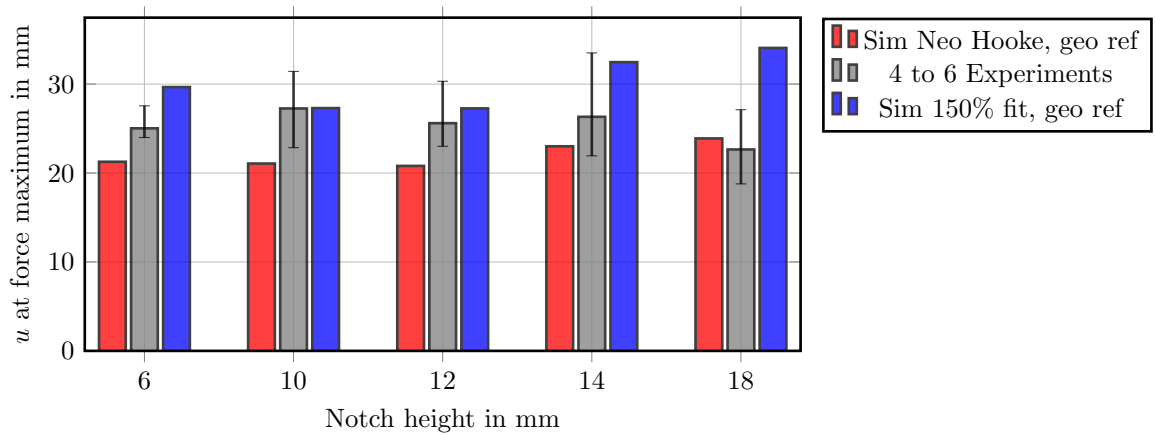

In order to investigate the crack path behavior, punctured strips are elongated with $200 \mathrm{mm/min}$ (related to traverse velocity) in a ZWICK Universal test machine till total failure. 
An inhomogeneous forming strain and stress field during the experiment is introduced by a circular hole of $8\mathrm{mm}$ diameter in the upper right part of the strip, cf. Figure \ref{geo_sim}. A similar experimental setup was proposed and investigated in \cite{Ozelo2012}.
Variations of the experiments are realized by varying the position of the given notch on the left side with a length of $1\mathrm{mm}$ at $6$, $10$, $12$, $14$ and $18 \mathrm{mm}$ height (from the bottom boundary). 
In all experiments the hole has a high influence on the crack path. For an initial position from $6$ to $10 \mathrm{mm}$ the crack path is diverted towards the hole. While for $6$ and $10 \mathrm{mm}$ the crack propagates below the hole towards the right edge, for $12$, $14$ and $18 \mathrm{mm}$ the crack is stopped for a short time by the hole, propagating afterwards nearly at the middle right inner edge of the hole towards the boundary edge of the specimen, see Figures \ref{fig:NHK00169_v2} and \ref{fig:crackpaths}. 
Comparing the start of the propagating fracture, the notch height only has minor influence on the force response. Although the distance between the initial notch and the circular hole is the shortest for the $14 \mathrm{mm}$ test, the resistance is quite strong. A clear trend between initial notch position and the start of the force drop in terms of the global displacement is not observable, cf. Figure \ref{fig:crackstart}. The high deviation in the (starting) crack behaviour between same experiments is typical for carbon black filled rubbers. Despite the deviations in the force drop, the crack path is quite stable within the group of same experiments for one notch height, cf. Figure \ref{fig:crackpaths}.\\

\begin{figure}[htbp!]
	\centering
		\includegraphics{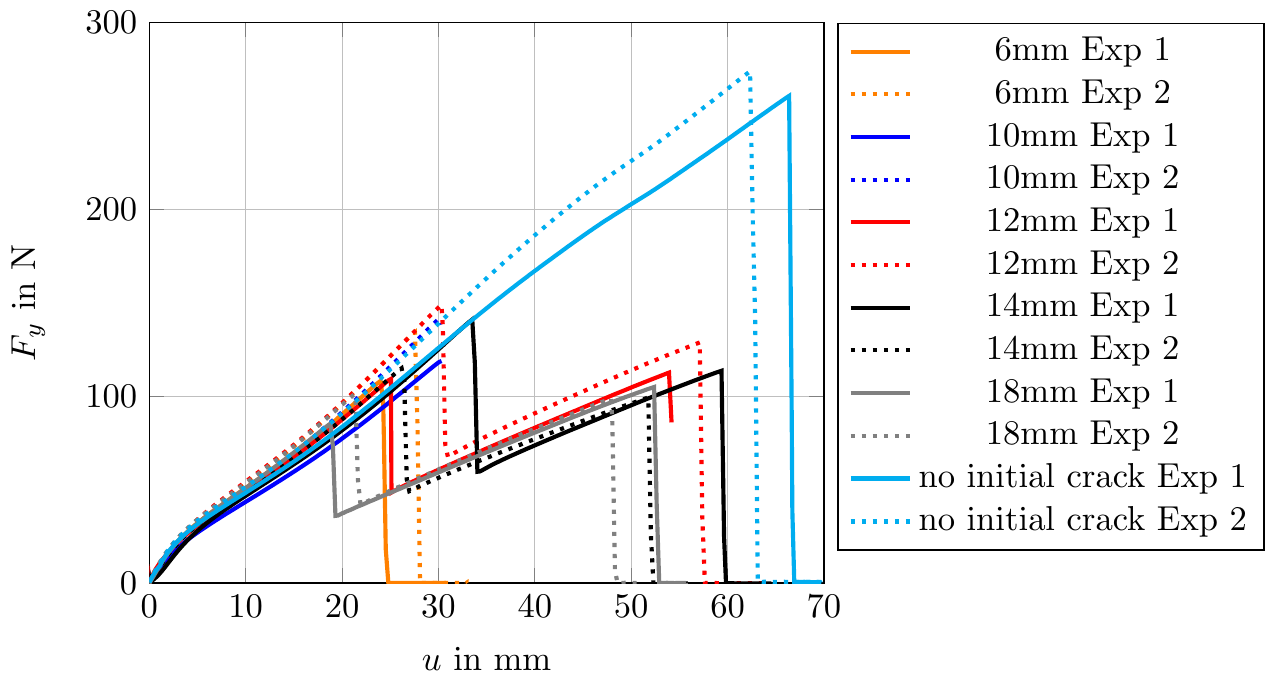}
	\caption{Two force-displacement curves for each notch height and test runs without an initial notch. Force response $F_y$ in relation to the traverse displacement $u$ measured in $y$-direction on the top boundary of the EPDM strips.}
	\label{fig:crackpathexperiments_fd}
\end{figure}

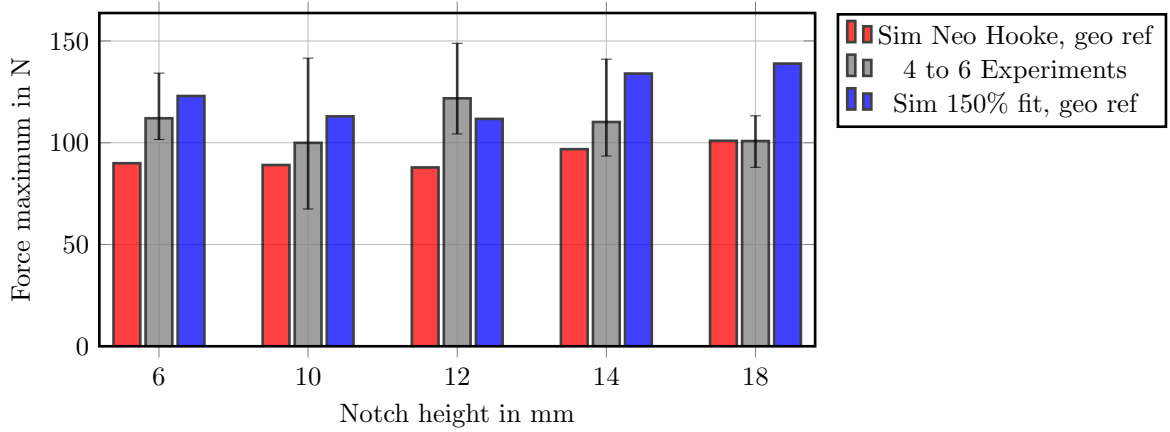
\begin{figure}[htbp!]
\centering
\begin{tikzpicture}
\begin{axis}[
    width=11cm,
    height=6cm,
    ylabel = Force maximum in $\mathrm{N}$,
    xlabel = Notch height in $\mathrm{mm}$,
    grid =major,
    ybar,
    symbolic x coords={6,10,12,14,18},xtick={6,10,12,14,18},
    ymin=0,
    legend pos=outer north east,
]
\addplot[fill=red,opacity=0.75] coordinates
{(6,89.9312) (10,89.0628) (12,87.8583) (14,96.8353) (18,100.96)};
\addlegendentry{Sim Neo Hooke, geo ref};
\addplot[fill=gray,opacity=0.75,error bars/.cd, y dir=both, y explicit] coordinates
{(6,112.011) -=(0,10.428) +=(0,22.178) (10,99.93576875) -=(0,32.52176875) +=(0,41.60823125) (12,121.8145) -=(0,17.5185) +=(0,27.0305) (14,110.18026) -=(0,16.71236) +=(0,30.89874) (18,100.8354) -=(0,12.9198) +=(0,12.4056)};
\addlegendentry{$4$ to $6$ Experiments};
\addplot[fill=blue,opacity=0.75] coordinates
{(6,122.955) (10,112.956) (12,111.691) (14,133.98) (18,138.876)};
\addlegendentry{Sim $150\%$ fit, geo ref};
\end{axis}
\end{tikzpicture}
\caption{Maximal loading force $F_y$ at the crack start measured on the top boundary, of experiments and finite element simulations. In gray (middle bar for each test), the averaged force maximum of $4$ to $6$ experiments in dependence of traverse displacement. In red and blue (left and right for each notch height), the numerically achieved force maxima are given.}
	\label{fig:crackmaxforce}
\end{figure}

In Figure \ref{fig:crackmaxforce}, the maximal loading forces at the crack start are given for the experiments on punctured strips with different notch height compared to the numerical results for each test. Concerning the experiments one can see a fully evolved crack to the opposite edge ($6$ and $10 \mathrm{mm}$) respectively or an evolved crack to the hole ($12, 14$ and $18 \mathrm{mm}$), cf. Figures \ref{fig:NHK00169_v2}, \ref{fig:crackpaths} and \ref{fig:crackpathexperiments_fd}. In red (left bar) and blue (right bar) in Figure \ref{fig:crackmaxforce}, the numerical results for the force maximum are plotted, based on the two coloured parameter settings from Table \ref{Tab:Youngsmodulus} for the notch heights of $6$ to $18\mathrm{mm}$.

In Figure \ref{fig:crackstart}, the experimental and numerical results of the traverse displacement at the first force maximum from Figure \ref{fig:crackmaxforce} are presented. The grey bars in the middle give the experimental weighted results including natural scattering of the experiments. In red (left bar) and blue (right bar), the results for the maximal displacement at the force maximum are plotted based on the parameter settings from Table \ref{Tab:Youngsmodulus}.


\section{A quasi-static phase-field model for nearly incompressible solids}\label{NewModel}
In this section, a quasi-static phase-field fracture model is proposed which is used to simulate crack propagation of the experiments described in the previous Section \ref{Experi}. In the following, the used mixed phase-field fracture model is given as well as a description on the numerical solving and the implementation of the finite element simulation, see Section \ref{NumCode}. In Section \ref{Simulation}, the numerical results concerning the punctured EPDM strips are presented and compared to the experimental results considering the crack paths and force-displacement curves.


\subsection{Notation, spaces and preliminaries}\label{notation}
In the following, we emanate from a two-dimensional, open and smooth domain $\Omega\subset \mathbb{R}^2$ where $\Gamma_D$ denotes the Dirichlet boundary. 
Let $\mathcal{I}$ be a loading interval $[0,T]$, where $T>0$ is the end time value.
A displacement function $u$ is given as $u:(\Omega \times \mathcal{I}) \to \mathbb{R}^2$. Further, a smooth indicator function $\varphi:(\Omega \times \mathcal{I}) \to [0,1]$ is named phase-field function with $\varphi=0$ in the broken, and $\varphi=1$ in the unbroken material.
The physics of the underlying problem ask to enforce crack irreversibility, i.e., that $\varphi$ is monotone non-increasing with respect
to $t\in \mathcal{I}$.\\
By $(a,b) := \int_\Omega a \cdot b\ dx$ for vectors $a, b$, the $L^2$ scalar-product is denoted. 
The Frobenius scalar product of two matrices of the same dimension is defined as $(A:B):= \sum_i \sum_j a_{ij} b_{ij}$ and therewith the $L^2$-scalar product is given by $(A,B) := \int_\Omega A : B\ dx$. 

To allow for a weak problem formulation, we consider a subdivision $0 =t_0 < t_1 < \ldots < t_N = T$ of the interval
$\mathcal{I}$. In each time step, we define approximations $(u^n,\varphi^n) \approx (u(t_n),\varphi(t_n))$ and hence the irreversibility condition is approximated by
$\varphi^{n}\leq\varphi^{n-1}$ for all $n = 1, \ldots, N$.
The phase-field space is $\mathcal{W} := H^1(\Omega)$ with a convex subset
$\mathcal{K}:=\{\psi\in \mathcal{W}\mid \psi\leq \varphi^{n-1}\leq 1\}$.
Further, we define the function spaces $\mathcal{V}:= (H_0^1(\Omega))^2:= \{w \in (H^1(\Omega))^2\mid w=0\mbox{ a.e. on }\Gamma_D\}$ and $\mathcal{U}:=L_2(\Omega)$.\\
In the following, the critical energy release rate is denoted by $G_c$.
To guarantee well-conditioning of the system of equations, a degradation function is defined as $g(\varphi):=(1-\kappa)\varphi^2 + \kappa,$
with a small regularization parameter $\kappa > 0$. 
The stress tensor $\sigma(u)$ is given by $\sigma(u) := 2 \mu E_{\text{lin}}(u) + \lambda \text{tr} (E_{\text{lin}}(u)) \textbf{I}$ with the Lam\'e coefficients $\mu,\lambda > 0$.
The linearized strain tensor therein is defined as $E_{\text{lin}}(u):=\frac{1}{2} (\nabla u + \nabla u^T)$. 
By $\textbf{I}$, the two-dimensional identity matrix is denoted.\\
We follow Wu \cite{wu2017unified} for a unified phase-field fracture model with the energy functional
\begin{align}
\begin{aligned}
E_{\epsilon}(u,\varphi) = \int_{\Omega} \frac{g(\varphi)}{2} \sigma(u): E_{\text{lin}}(u)\, dx
+ \int_{\Omega} \frac{G_c}{\pi} \frac{2(1-\varphi) - (1-\varphi)^2}{\epsilon}\, dx + \int_{\Omega} \frac{G_c}{\pi} \epsilon |\nabla \varphi|^2\, dx,\label{functional}
\end{aligned}
\end{align}
where $\epsilon >0$ is the usual phase-field regularization parameter, describing the bandwidth of the transition zone between broken and unbroken material. 

\begin{remark}[Wu's energy functional]
Ambrosio and Tortorelli \cite{ambrosio1990approximation,ambrosio1992approximation} proposed an elegant way to approximate the Mumford-Shah functional \cite{mumford1989optimal} for image segmentation by an elliptic functional defined on Sobolev spaces. The key idea was to replace the sharp lower-dimensional crack by a smoothed indicator function while $\Gamma$-convergence to the Mumford-Shah functional is guaranteed for a length scale $\epsilon\to 0$. This in turn satisfies Griffith's law of crack propagation.
Simultaneously, Bourdin et al. \cite{bourdin2000numerical} proposed a regularized energy functional for brittle fracture. The so called AT$_2$ functional (named as in \cite{tanne2018crack}) is defined as
\begin{align}
\begin{aligned}
{\text{AT}_2}:\quad E_{\epsilon}(u,\varphi) = \int_{\Omega} \frac{g(\varphi)}{2} \sigma(u): E_{\text{lin}}(u)\, dx
+ \int_{\Omega} \frac{G_c}{2} \frac{(1-\varphi)^2}{\epsilon}\, dx + \int_{\Omega} \frac{G_c}{2} \epsilon |\nabla \varphi|^2\, dx.\label{at2}
\end{aligned}
\end{align}
Later in 2014, Bourdin et al. \cite{bourdin2014morphogenesis} introduced a very similar functional with a stress-softening behavior and where the damage model remains purely elastic without damage until the stress
reaches the critical value. This functional for complex cracks is defined as
\begin{align}
\begin{aligned}
{\text{AT}_1}:\quad E_{\epsilon}(u,\varphi) = \int_{\Omega} \frac{g(\varphi)}{2} \sigma(u): E_{\text{lin}}(u)\, dx
+ \int_{\Omega} \frac{3G_c}{8} \frac{1-\varphi}{\epsilon}\, dx + \int_{\Omega} \frac{3G_c}{8} \epsilon |\nabla \varphi|^2\, dx.\label{at1}
\end{aligned}
\end{align}
Aside from different coefficients, Wu's energy functional in Equation \eqref{functional} in comparison to the functionals in the Equations \eqref{at2} and \eqref{at1} uses a combination of a linear and a quadratic part in the second crack energy term, which yields the theoretically and numerically useful property of a finite support for a localized phase-field \cite{wu2017unified}. The numerical consequences of this choice are described in Section \ref{comment}.
\end{remark}

The Euler-Lagrange equations of Wu's energy functional in Equation \eqref{functional} look like follows:
\begin{problem}[Coupled variational inequality system of Wu's energy functional]\label{prob_1}
Find $\varphi \in \mathcal{K}$ and $u \in \mathcal{V}$, such that 
\setlength\jot{0.4cm}
 \begin{align}
 \begin{aligned}
 \left(g(\varphi) 2\mu E(u),E(w)\right)
  + \left(g(\varphi) \lambda \text{tr} E_{\text{lin}}(u) \textbf{I}, E_{\text{lin}}(w)\right)=&\ 0\quad \forall w \in \mathcal{V},\\
  (1- \kappa) (\varphi 2\mu E_{\text{lin}}(u):E_{\text{lin}}(u),\psi-\varphi) 
+ (1-\kappa)(\varphi \lambda \nabla \cdot u\ \textbf{I}:E_{\text{lin}}(u),\psi-\varphi)\\
+ \frac{2 G_c}{\pi} (-\frac{1}{\epsilon} \varphi,\psi-\varphi)
+ \frac{2 G_c}{\pi} \epsilon(\nabla \varphi,\nabla(\psi-\varphi))\geq&\ 0\quad \forall \psi \in \mathcal{K}. \label{wu_pde}
 \end{aligned}
 \end{align}
 \end{problem}
Based on the coupled variational inequality system \eqref{wu_pde}, a mixed phase-field fracture model is developed in the following.
\subsection{Mixed phase-field fracture model}\label{mixedModel}
For a mixed form of the problem, we define a hydro-static pressure $p\in \mathcal{U}$
\begin{align*}
 p:= \lambda \text{tr} \left(E_{\text{lin}}(u)\right),
\end{align*}
such that the pure elasticity equation in weak form reads as:
\begin{problem}[Pure elasticity]\label{prob_2}
Assume $\varphi \in \mathcal{K}$ to be given. Find $u \in \mathcal{V}$ and $p \in \mathcal{U}$ such that 
 \begin{align}
 \quad\quad\quad\quad\quad\quad\quad\quad g(\varphi) 2\mu \left(E_{\text{lin}}(u),E_{\text{lin}}(w)\right)+g(\varphi) (\nabla \cdot w, p)=&\ 0\quad \forall w \in\ \mathcal{V},\\
g(\varphi) (\nabla \cdot u,q) - \frac{1}{\lambda} (p,q)=&\ 0\quad \forall q\in\ \mathcal{U}.\label{incomp}
 \end{align}
 \end{problem}
 
 We refer the reader to \cite{mang2020phase} for details on the inf-sup stable mixed form of the elasticity equation.
 
\begin{remark}[Degradation in the crack]
  To avoid non-physical pressure values in the inner fracture zone, the degradation term $g(\varphi)$ is neglected in the second term of Equation \eqref{incomp}. It follows, if we are in the broken zone $\varphi=0$, we are not fully divergence-free. This modification is an extension of our own recent model \cite{mang2020phase}.
\end{remark}

The variational formulation in Problem \ref{prob_1} and Problem \ref{prob_2} of the isotropic model given by Bourdin et al. \cite{bourdin2000numerical} does not allow to distinguish between fracture behaviour in tension and compression \cite{ambati2015review}. In the following, we consider the volumetric and deviatoric contributions of the elastic energy density separated into $\sigma^+$ and $\sigma^-$ to prevent interpenetration of the crack faces under compression~\cite{ambati2015review}, along to Amor et al.~\cite{amor2009regularized}. For this reason, the positive part of the pressure $p^+ \in L_2(\Omega)$ as well as the positive part of $E_{\text{lin}}^+(u)$ have to be defined as
$p^+ := \max \{p,0\}$ and $E_{\text{lin}}^+(u) := \max \{E_{\text{lin}}(u),0\}$ as the maximum function,
such that the stress tensor $\sigma(u,p)$ is split into:
\begin{align*}
\begin{aligned}
 \sigma^+(u,p):=&\ \mu \max\left\{0, \text{tr} \left(E_{\text{lin}}^+(u)\right)\right\} \textbf{I} 
 + 2\mu \left( E_{\text{lin}}^+(u) - \frac{1}{3} \text{tr} \left(E_{\text{lin}}^+(u)\right)\textbf{I}\right) + p^+ \textbf{I},\\
 \sigma^-(u,p):=&\ \mu \left(\text{tr} \left(E_{\text{lin}}^+(u)\right) - \max\left\{0, \text{tr} \left(E_{\text{lin}}^+(u)\right)\right\}\right) \textbf{I} + (p - p^+) \textbf{I}.
 \end{aligned}
\end{align*}
\begin{remark}
The mixed quasi-static phase-field fracture model developed and used in previous contributions \cite{mang2020phase,basava2020adaptive,mang2020adaptive}, is derived from the anisotropic phase-field fracture model with a spectral decomposition of the strain tensor along to Miehe et al. \cite{MieWelHof10a}. Based on this model promising numerical results considering benchmark tests in compressible and (nearly) incompressible solids could be presented. In the work on hand, the phase-field fracture model is applied to simulate a more complex crack propagation in nearly incompressible EPDM rubber, where the volumetric-deviatoric split of Amor et al. \cite{amor2009regularized} allows more realistic numerical results. We refer the reader to Ambati et al. \cite{ambati2015review}, who discuss the very similar numerical results for both energy splitting approaches. Nevertheless, the computational effort of the volumetric-deviatoric strain energy split is minor than for the spectral decomposition of the strain tensor.
\end{remark}

From that, our newly introduced mixed phase-field fracture problem based on Wu's energy functional and strain energy splitting of Amor et al. in incremental form is formulated as:

\begin{problem}[Mixed phase-field formulation]\label{formMixed}
Let the previous loading step data $\varphi^{n-1} \in \mathcal{K}$ be given.
Find $u:= u^n \in \{u_D+\mathcal{V}\}$, $p:=p^n \in \mathcal{U}$ and $\varphi:=
\varphi^n \in \mathcal{K}$ for the loading steps
$n=1,2,\ldots, N$ such that 
\setlength\jot{0.3cm}
 \begin{align}
 g(\varphi^{n-1})\left(\sigma^+(u,p),E_{\text{lin}}(w) \right) + \left(\sigma^-(u,p),E_{\text{lin}}(w)\right)=&\ 0\quad \forall\ w\in \mathcal{V},\\ \label{elast_1}
  g(\varphi^{n-1})\left(\nabla \cdot u,q\right) - \frac{1}{\lambda} (p,q)=&\ 0\quad \forall\ q\in \mathcal{U},\\ \label{elast_2}
  (1- \kappa)\left(\varphi \sigma^+(u,p) : E_{\text{lin}}(u),\psi-\varphi\right) + \frac{2G_c}{\pi} \left(-\frac{1}{\epsilon} \varphi,\psi- \varphi\right)\nonumber \\
  + \frac{2G_c}{\pi} \epsilon \left(\nabla \varphi,\nabla(\psi-\varphi) \right)\geq&\ 0\quad \forall\ \psi \in \mathcal{K}.
 \end{align}
\end{problem}
\begin{remark}
In the elasticity Equations \eqref{elast_1} and \eqref{elast_2}, time-lagging is used in the phase-field variable $\varphi$ to obtain a convex functional. 
Specifically, we notice that of course a time-discretization error of 
the order of the time step size (loading increment) arises. For fast growing 
cracks this makes a significant difference as shown in \cite{wick2017error}. 
By applying subiterations, the discretization error can be reduced \cite{wick2020multiphysics}[Chapter 7.7].
\end{remark}
To the best of our knowledge, the phase-field fracture model derived from Wu's energy functional in Equation \eqref{functional} with a mixed variational formulation of the elasticity equation, has not been considered in the literature so far.
The numerical handling of the coupled phase-field fracture system with an inequality constraint in Problem \ref{formMixed} with the help of a monolithic solving approach including the discretization and two different refinement strategies are discussed in the following.


\subsection{Numerical treatment, implementation, and programming code}\label{NumCode}
The crack irreversibility, i.e., the variational inequality in the phase-field part is fulfilled with a primal-dual active set method; see \cite{heister2015primal} for further details. The overall nonlinear problem is solved with a Newton method and derived as a joint quasi-monolithic approach combined with the primal-dual active set method \cite{heister2015primal}.
We notice that the nonlinear monolithic solution is challenging 
and different studies have appeared \cite{wick2017error,heister2015primal,GeLo16,wick2017modified,KoKr20,KRISTENSEN2020102446,WAMBACQ2021113612,WU2020102440,WU2020112704,JoLaWi20}.

For the spatial discretization, we employ a Galerkin finite element method in each incremental step, where the domain $\Omega$ is partitioned into quadrilaterals. 
To fulfill a discrete inf-sup condition, stable Taylor-Hood elements with 
bi-quadratic shape functions ($Q_2$) for the displacement field $u$ and bilinear shape functions ($Q_1$)
for the pressure variable $p$ and the phase-field variable $\varphi$ are used as in \cite{mang2020phase}.

The implementation of the proposed problem formulation is derived from the open-source code \url{https://github.com/tjhei/cracks}. Further details on the code are given in \cite{heister2018parallel,HeiWi20}. This 'pfm-cracks' project is based on deal.II~\cite{arndt2020deal} and provides to simulate crack propagation in elastic and porous media. Starting from the classical phase-field fracture model with two unknowns, the displacements and the phase-field variable discretized with $Q_1 Q_1$ elements, we implemented the mixed problem formulation with Taylor-Hood elements $Q_2 Q_1$ elements for the elasticity pair $(u,p)$. We end up with three unknowns $u,p$ and $\varphi$ discretized with $Q_2 Q_1 Q_1$ finite elements and solve the variational inequality system by help of a primal-dual active set scheme. Further, we set up the new test considering the punctured EPDM strips till total failure with proper boundary and initial conditions, and implemented the required functionals of interest to allow comparing numerical and experimental results.

For a spatial convergence study in Section \ref{space_conv}, a predictor-corrector scheme is used along to~\cite[Chapter 4]{heister2015primal} for adaptive mesh refinement. This refinement scheme allows to refine the mesh locally depending on a propagating fracture with a chosen threshold (here $\rho=0.7$) for the phase-field variable: when the crack phase-field on a coarse cell has values smaller than $0.7$, the mesh on this cell is refined and the incremental step is computed again on an adaptive mesh, where the cell with $\varphi<0.7$ is refined in dependence of its neighbour cells. This allows to give a sharp crack while having a reasonable amount of workload, because the mesh is just refined, if we have a change in the phase-field variable which pictures the crack area. In further numerical tests in Section \ref{comparison}, geometrically prerefined meshes are used to avoid observed mesh locking around the circular hole and to guarantee a good resolution also towards the top boundary where we evaluate the loading force as one quantity of interest.



\section{Numerical simulation of crack propagation in punctured EPDM strips}\label{Simulation}
Crack propagation in Section \ref{Experi} is observed via stretching the punctured EPDM strips until total failure with a speed of $200 \mathrm{mm/min}$.
For the numerical simulation of the crack propagation in punctured strips in a two-dimensional setup (plane stress assumption) with help of the phase-field fracture model from Section \ref{mixedModel},
we reduce the geometry to the area of interest between the bulges on the bottom and top part where the specimens are fixed. The geometry is given in Figure \ref{geo_sim}. Considering holes in a material combined with phase-field fracture modeling is still a challenging task and not fully understood; see \cite[Section 8.2]{wu2018phase} for further comments on the difficulty of inclusions.  Homogeneous Dirichlet boundary conditions $u_y=0$ are chosen on the top boundary and the strips are fixed in horizontal direction on the top and bottom boundary via $u_x=0$. The following boundary condition characterizes the loading force on the bottom boundary $\Gamma_{\text{force}}$ in $y$-direction:
 \begin{align}
 u_y = t \cdot 200 \mathrm{mm}/\mathrm{min},\ \text{for}\ t \in \mathcal{I}:= [0; \text{total failure}],\label{gamma_force}
\end{align}
where $t$ denotes the total time. In the quasi-static context, the time interval $\mathcal{I}$ is divided into incremental steps of size $\delta t$. Further, the phase-field is fixed via $\varphi=0$ in the given notch as an initial boundary condition. 
\begin{remark}[Numerical handling of given notch]
We decided to handle the given notch on the left side of size $1\mathrm{mm}$ with an initial condition $\varphi=0$ such that the material is already broken in the notch area. In addition, we describe the initial crack by doubling the
degrees of freedom on the respective faces similar to Wick \cite[Section 5.1.]{wick2017modified}. This allows the material to open in the notch and the maximal stress is obtained in the notch tip, which in turn imitates the observed opening of the notch in the running experiments.
\end{remark}
The numerical 
parameters are listed in Table \ref{num_params}.

\begin{figure}[htbp!]
\begin{minipage}{0.4\textwidth}
 \begin{tikzpicture}[xscale=1.5,yscale=1.5]
\draw[fill=gray!30] (0,0) -- (0,2.8) -- (2,2.8) -- (2,0) -- (0,0);
\draw (0.3,0) -- (0.15,-0.1);
\draw (0.5,0) -- (0.35,-0.1);
\draw (0.7,0) -- (0.55,-0.1);
\draw (0.9,0) -- (0.75,-0.1);
\draw (1.1,0) -- (0.95,-0.1);
\draw (1.3,0) -- (1.15,-0.1);
\draw (1.5,0) -- (1.35,-0.1);
\draw (1.7,0) -- (1.55,-0.1);
\draw (1.9,0) -- (1.75,-0.1);
\node[blue] at (1.0,3.1) {$\Gamma_{\text{top}}$};
\draw (0.3,2.9) -- (0.15,2.8);
\draw (0.5,2.9) -- (0.35,2.8);
\draw (0.7,2.9) -- (0.55,2.8);
\draw (0.9,2.9) -- (0.75,2.8);
\draw (1.1,2.9) -- (0.95,2.8);
\draw (1.3,2.9) -- (1.15,2.8);
\draw (1.5,2.9) -- (1.35,2.8);
\draw (1.7,2.9) -- (1.55,2.8);
\draw (1.9,2.9) -- (1.75,2.8);
\draw[fill=white] (1.2,1.8) circle[radius=0.4cm];
\node at (1.2,1.8) {\textbf{$\times$}};
\node at (1.2,1.65) {\scriptsize{$(12,18)$}};
\node at (1.2,2.35) {$\diameter 8 \mathrm{mm}$};
\draw[maroon] (0,1.0) -- (0.2,1.0);
\node[maroon] at (0.4,0.9) {notch};
\draw[<->] (-0.25,0) -- (-0.25,2.8);
\node at (-0.9,1.4) {$28 \mathrm{mm}$};
\draw[<->] (0,-0.25) -- (2.0,-0.25);
\node at (1.0,-0.5) {$20 \mathrm{mm}$};
\draw[blue, thick,->] (1.8,0.0) -- (1.8,-0.4);
\node[blue] at (1.9,-0.6) {$\Gamma_{\text{force}}$};
 \end{tikzpicture}
\captionof{figure}{Geometry and boundary conditions of punctured EPDM strips. Loading force on $\Gamma_{\text{force}}$ defined in Equation \eqref{gamma_force}.}\label{geo_sim}
\end{minipage}
\begin{minipage}{0.6\textwidth}
\begin{tabular}{|c|l|c|}\hline
 \rowcolor{gray!30} \multicolumn{1}{|c|}{Parameter} & \multicolumn{1}{|c|}{Description}  & \multicolumn{1}{|c|}{Value} \\ \hline \hline
    $h$  & discretization parameter (coarsest) & $0.3 \mathrm{mm}$ \\ \hline 
    $\#$ ref & number of adaptive refinement steps & $0$ to $3$ \\ \hline
    $\rho$ & phase-field threshold for refinement & $0.7$ \\ \hline
   $\epsilon$ & bandwidth  & $2h$\\ \hline
   $\delta t$ & incremental size & $10^{-2} \mathrm{s}$ \\ \hline
      $\kappa$ & regularization parameter & $0.01h$ \\ \hline
 \end{tabular}
\captionof{table}{Numerical parameters for the phase-field fracture simulation of punctured EPDM strips. Parameters '$\#$ ref' and '$\rho$' required for the adaptive test runs in Section \ref{space_conv}. The numerical parameters $h, \epsilon, \delta t$ and $\kappa$ are valid for all numerical test runs.}\label{num_params}
\end{minipage}
\end{figure}

The setting of the numerical parameters in Table \ref{num_params} is valid for all numerical results considering the spatial convergence study with the $6\mathrm{mm}$ test using adaptive refinement. The incremental size $\delta t$ and the size of $\kappa$ is the same for all numerical tests. Related to the material parameter setting, we choose two cases of Figure \ref{fig:EModul} and Table \ref{Tab:Youngsmodulus}. The coloured material parameters in Table \ref{Tab:Youngsmodulus} reflect the parameter settings based on two different assumptions for a linear $E$: the Young's modulus estimated via a Neo Hookean model and a fit of $150\%$ of the strain. In all following numerical results and especially in the figures of this section the setting of the material parameters is differed via naming it '$150\%$ fit' or 'Neo Hooke', respectively.
The bulk modulus $K = 2595 \mathrm{MPa}$ listed in Table \ref{Tab:Youngsmodulus} is specified via an averaged value of the fitted values in Figure \ref{fig:compression}. The critical energy release rate is adopted from Figure \ref{fig:tearingenergyspecimen} such that $G_c=17.0\mathrm{N/mm}$.\\

\textbf{Quantities of interest.}
To compare the experimental and numerical results, we consider two quantities of interest: first, the solution of the phase-field after total failure and the crack paths in the punctured strips. As a second quantity, we compare for proper test cases the force-displacement curves of the experiments and simulations to retrace the crack propagation process.

Relative to the force-displacement curves, the load vector on the top boundary $\Gamma_{\text{top}}$ is evaluated via 
\begin{align}
 (F_x,F_y):=\frac{1}{|\Gamma_{\text{top}}|}\int_{\Gamma_{\text{top}}} \sigma(u) n\ ds,\label{loading}
\end{align}
with the length $|\Gamma_{\text{top}}|=20\mathrm{mm}$ of the top boundary, the stress tensor $\sigma(u) := 2 \mu E_{\text{lin}}(u) + \lambda \text{tr}  (E_{\text{lin}}(u)) \textbf{I}$ and the normal vector $n$. 
Within the EPDM tests we are interested in the loading force $F_y$ on $\Gamma_{\text{top}}$, where the force response is also measured in the experiments.\\

\subsection{Spatial convergence study}\label{space_conv}
\begin{figure}[htbp!]
\centering
\includegraphics[width=0.9\textwidth]{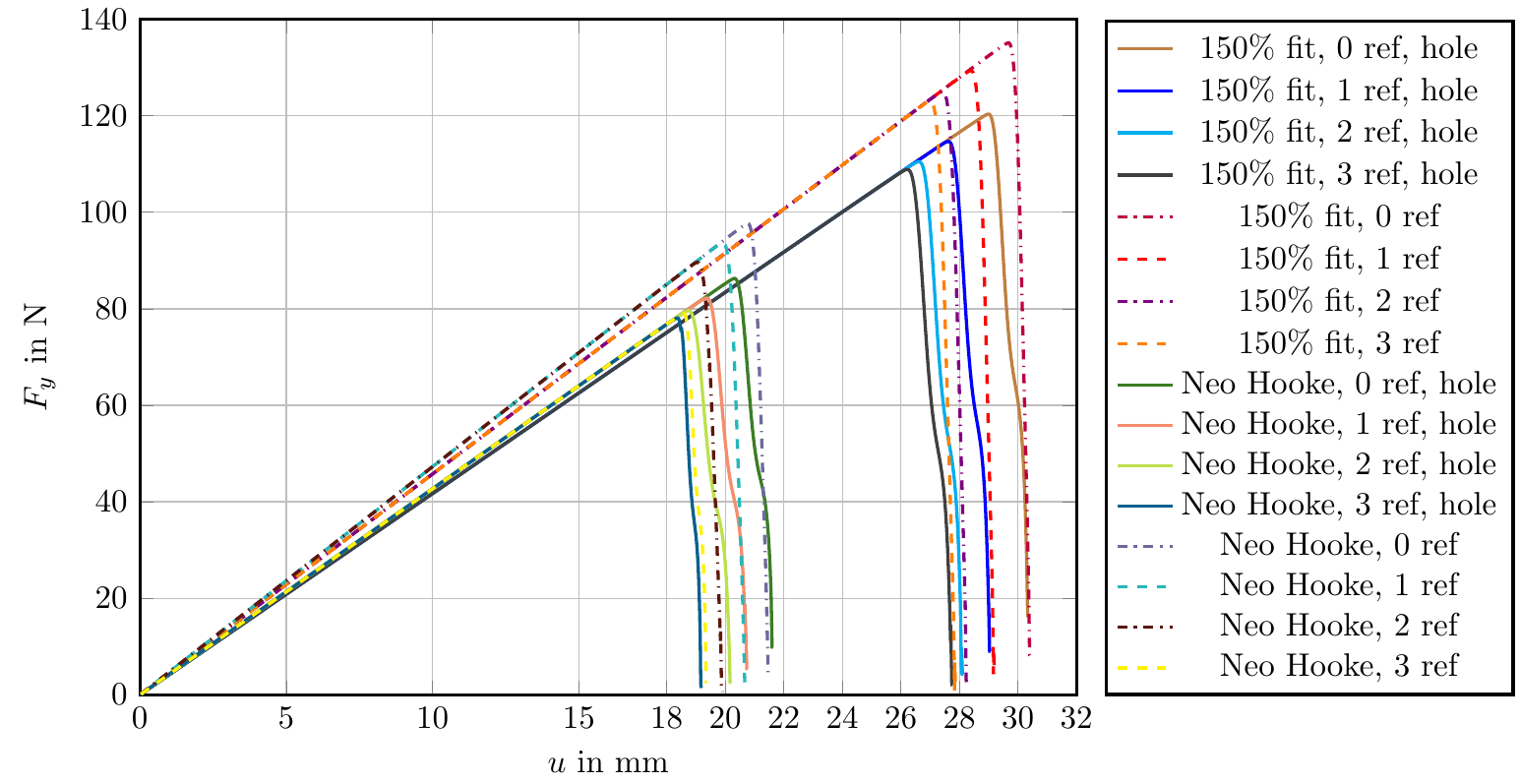}
\caption{Numerical refinement study: force-displacement curves (displacement $u$ in $y$-direction versus force $F_y$) for the EPDM benchmark test based on the material parameters of Table \ref{Tab:Youngsmodulus} ($150\%$ fit of strain) compared the material parameters of Table \ref{Tab:Youngsmodulus} with a notch at $6 \mathrm{mm}$, $h_{\text{start}}=0.3$ and $\epsilon =2h_{\text{start}} = 0.6 \mathrm{mm}$ fixed. $G_c = 17.0\mathrm{N/mm}$. The computation of the force response $F_y$ on the top boundary is defined in Equation \eqref{loading}.}\label{compare_hole_no_hole_150fit_NeoHooke}
\end{figure}

\begin{figure}[htbp!]
\centering
\includegraphics[height=3.5cm]{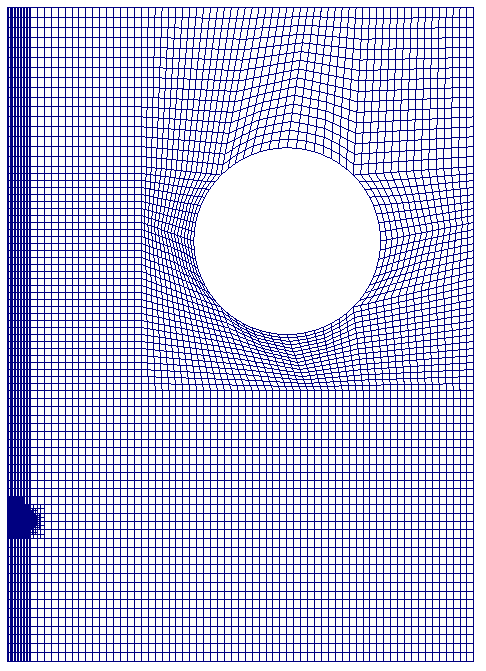}
\hfill
\includegraphics[height=3.5cm]{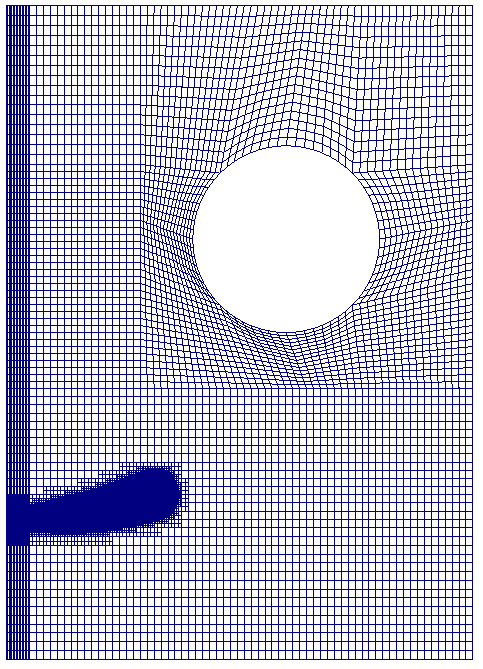}
\hfill
\includegraphics[height=3.5cm]{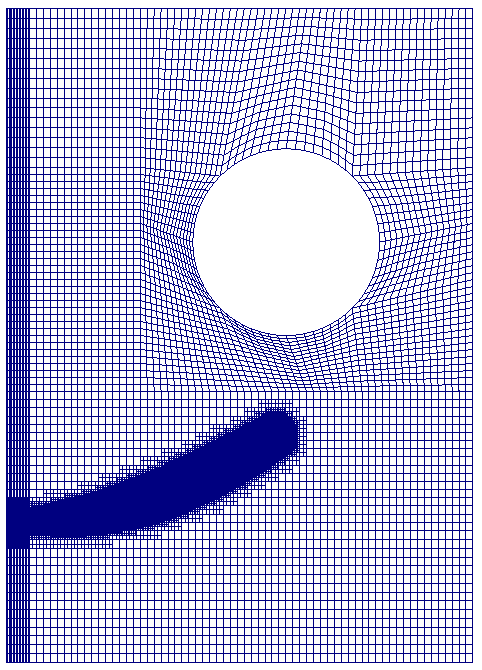}
\hfill
\includegraphics[height=3.5cm]{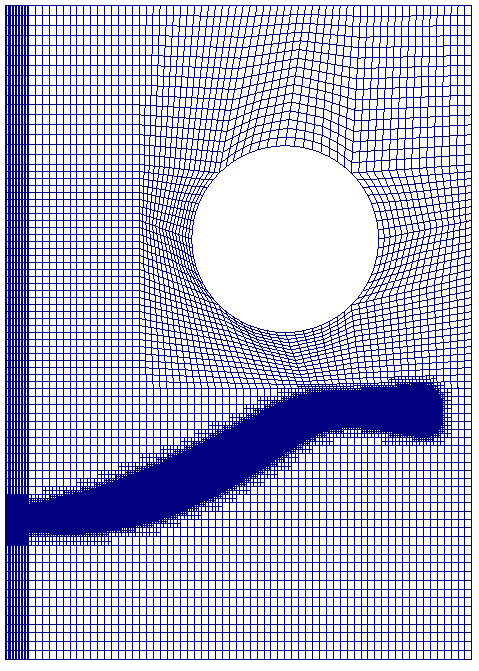}
\hfill
\includegraphics[height=3.5cm]{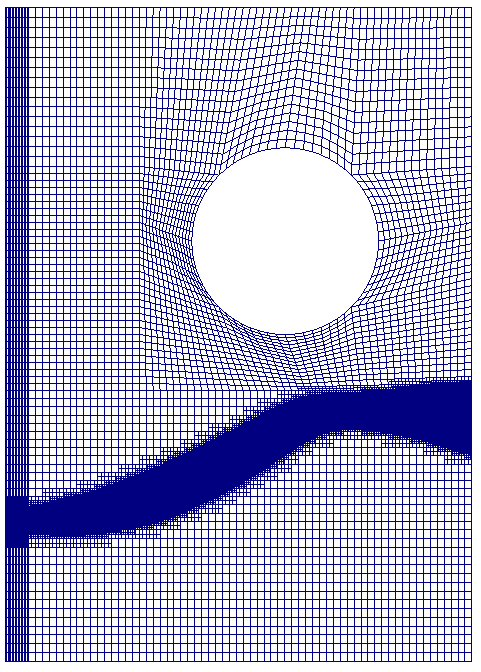}
\caption{Adaptively refined meshes at five certain time points for the numerical simulation based on the Neo Hookean material parameter setting from Table \ref{Tab:Youngsmodulus} and three steps of adaptive refinement steps via the predictor-corrector scheme from \cite{heister2015primal}, cf. the blue-green curve in Figure \ref{compare_hole_no_hole_150fit_NeoHooke}. In the last snapshot, the mesh contains $79,954$ DoF for the solid displacements and $10,160$ DoF for the pressure and the phase-field variable. }\label{adapt_mesh}
\end{figure}

To study spatial convergence of the finite element simulation, we start with a series of tests for one of the five test setups: the given notch is at a height of $6 \mathrm{mm}$, where experimentally we do not expect the crack propagating towards the hole but in a 'S'-curve from the left to the right. For this example, four test runs are conducted, first on a global prerefined mesh with a discretization parameter $h_{\text{start}}=0.3\mathrm{mm}$ and three tests with $1$ to $3$ adaptive refinement steps. The numerical results are given in Figure \ref{space_conv}. The number of degrees of freedom (DoF) for all conducted tests runs of spatial convergence are listed in Table \ref{dofs}.  While refining the area around the crack adaptively along to the predictor-corrector scheme from \cite{heister2015primal}, or in other words, while decreasing $h$ locally, the length scale parameter $\epsilon$ is fixed. To get an idea how the predictor-corrector scheme works within the propagating fracture, for one test of the series ($6\mathrm{mm}$, Neo Hooke, hole), the meshes at five certain loading points are given in Figure \ref{adapt_mesh}. Consequently, for all computations considered in Figure \ref{compare_hole_no_hole_150fit_NeoHooke}, it holds $\epsilon = 0.6\mathrm{mm}$. Further, the results for both parameter settings '$150\%$ fit' and 'Neo Hooke' are given.\\

Beside of spatial convergence, we want to avoid locking effects that could arise around the circular hole. The test with a notch at $6 \mathrm{mm}$ height gives a crack path away from the hole, even if the inclusion has an impact on the crack pattern. The refinement series is also conducted for both different material parameter settings (Neo Hooke and $150\%$ strain fit) without the inclusion.

In total, three points are observed in Figure \ref{compare_hole_no_hole_150fit_NeoHooke}: 
\begin{itemize}
    \item We compare the results based on the Neo Hookean material parameters from Table \ref{Tab:Youngsmodulus} with the computation based on the $150\%$ strain fit assumption.
    \item The results are given based on the geometry in Figure \ref{geo_sim} with a given notch and a circular hole in the upper right part of the strips compared to results based on the same geometry with a given notch but without a hole. The results in Figure \ref{compare_hole_no_hole_150fit_NeoHooke} for the punctured strips are marked with the attribute 'hole' in the legend.
    \item The third point discussed in Figure \ref{compare_hole_no_hole_150fit_NeoHooke}, is spatial convergence of the force-displacement curves with an increasing number of adaptive refinement steps, while the crack width $\epsilon$ is fixed as in Table \ref{num_params} even if $h$ is getting smaller in the crack area. In the legend in Figure \ref{compare_hole_no_hole_150fit_NeoHooke} the number of refinement steps ('$\#$ ref') is given for each force-displacement curve.
\end{itemize}
It can be observed that the gradient towards the force maximum in the results based on the punctured strips is steeper than without a circular inclusion. From this, one could follow, that the inclusion makes the material more elastic and decreases the material stiffness. Further, a significant difference between the force-displacement curves of the tests based on the $150\%$ fit or Neo Hookean material parameters can be observed. All Neo Hookean tests have a lower maximal force, which leads to a smaller displacement when the crack starts propagating, cf. Figures \ref{fig:crackmaxforce} and \ref{fig:crackstart}. Further in Figure \ref{compare_hole_no_hole_150fit_NeoHooke} one can observe spatial convergence within a series of tests from $0$ to $3$ adaptive refinement steps.

\begin{figure}[htbp!]
\centering
\begin{tabular}{|c|r|r|}\hline
 \rowcolor{gray!30} \multicolumn{1}{|c|}{Test (name)} & \multicolumn{1}{|c|}{$\#$DoF $u$}  & \multicolumn{1}{|c|}{($\#$DoF $p$) $=$ ($\#$DoF $\varphi$) } \\ \hline \hline
    $6\mathrm{mm}$, $150\%$ fit $\&$ Neo Hooke, 0 ref, hole  & 49,984 & 6,352 \\ \hline
    $6\mathrm{mm}$, $150\%$ fit, 1 ref, hole & 65,284 & 8,277 \\ \hline
    $6\mathrm{mm}$, $150\%$ fit, 2 ref, hole & 117,400 & 1,4815 \\ \hline
    $6\mathrm{mm}$, $150\%$ fit, 3 ref, hole & 304,786 & 38,288 \\ \hline\hline
    $6\mathrm{mm}$, $150\%$ fit, 0 ref  & 51,906 & 6,577 \\ \hline
    $6\mathrm{mm}$, $150\%$ fit, 1 ref & 70,376 & 8,898 \\ \hline
    $6\mathrm{mm}$, $150\%$ fit, 2 ref & 135,220 & 17,028 \\ \hline
    $6\mathrm{mm}$, $150\%$ fit, 3 ref & 369,498 & 46,359 \\ \hline\hline
    $6\mathrm{mm}$, Neo Hooke, 1 ref, hole & 65,420 & 8,294 \\ \hline
    $6\mathrm{mm}$, Neo Hooke, 2 ref, hole & 119,396 & 15,067 \\ \hline
    $6\mathrm{mm}$, Neo Hooke, 3 ref, hole (Fig. \ref{adapt_mesh}) & 293,572 & 36,880 \\ \hline\hline
    $6\mathrm{mm}$, Neo Hooke, 0 ref  & 51906 & 6,577 \\ \hline
    $6\mathrm{mm}$, Neo Hooke, 1 ref & 70,446 & 8,907 \\ \hline
    $6\mathrm{mm}$, Neo Hooke, 2 ref & 136,742 & 17,308 \\ \hline 
    $6\mathrm{mm}$, Neo Hooke, 3 ref & 371,840 & 46,650 \\ \hline\hline
    Figure \ref{18mm_compare}  & 237,436 & 29,877 \\ \hline
    Figure \ref{14mm_compare} & 391,656 & 49,192 \\ \hline
    Figure \ref{12mm_compare}  & 427,932 & 53,731 \\ \hline
    Figure \ref{10mm_compare}  & 437,544 & 54,934 \\ \hline
    Figure \ref{6mm_compare}  & 581,724 & 72,979 \\ \hline
 \end{tabular}
\captionof{table}{Number of degrees of freedom (DoF) for all numerical test runs in Figure \ref{compare_hole_no_hole_150fit_NeoHooke}, and the Figures \ref{18mm_compare} to \ref{6mm_compare} for the displacement $u$ ($Q_2$-elements) and the pressure variable $p$ ($Q_1$-elements), which has the same number of DoF as the phase-field variable $\varphi$ ($Q_1$-elements). The number of degrees of freedom for the test with adaptive predictor corrector mesh refinement from Figure \ref{compare_hole_no_hole_150fit_NeoHooke} are given for the mesh at total failure in the last computed incremental step.}\label{dofs}
\end{figure}

\begin{figure}[htbp!]
\centering
\includegraphics[width=0.9\textwidth]{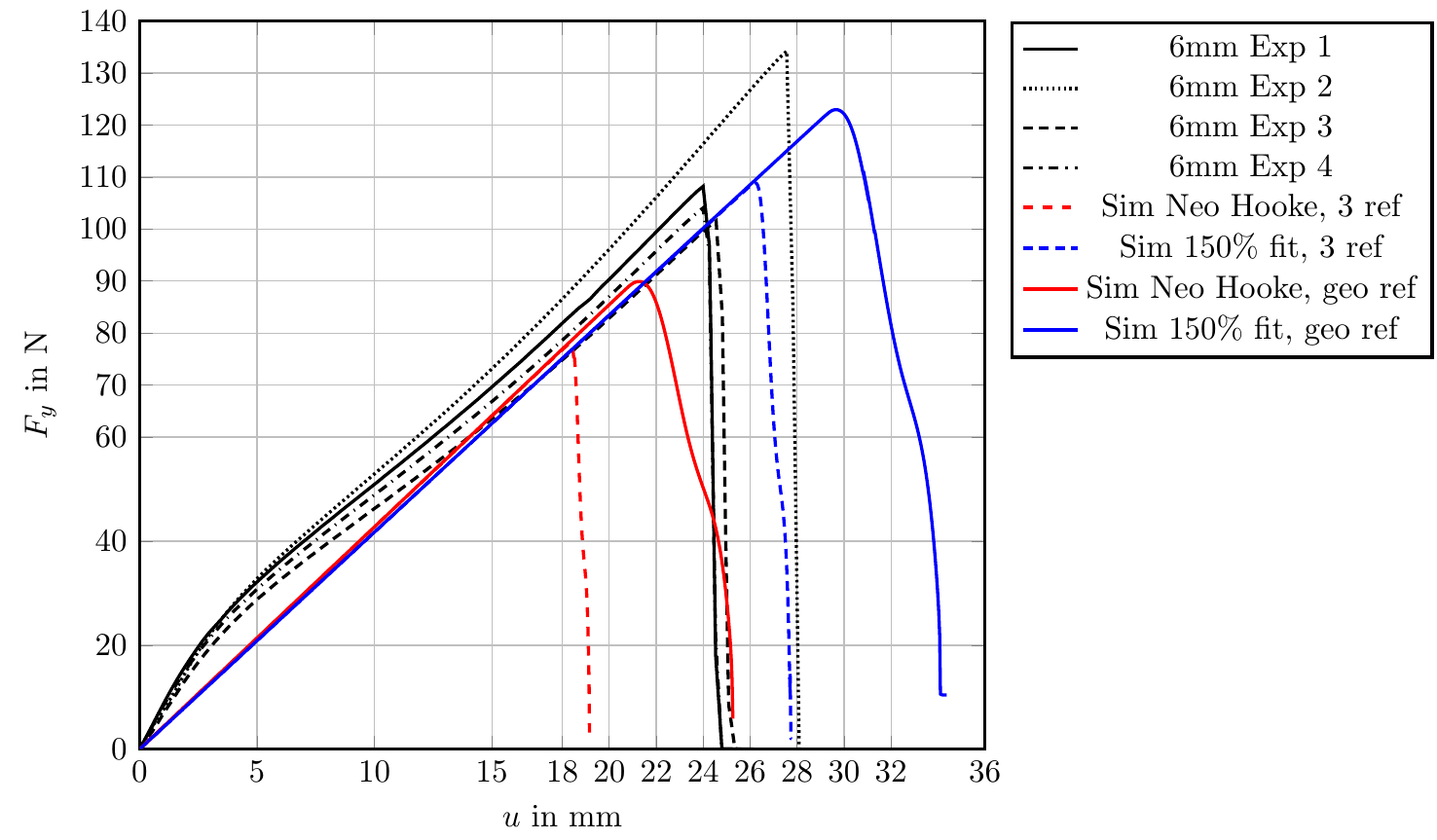}
\caption{Comparison of the force-displacement curves (displacement $u$ in $y$-direction versus force $F_y$) for the samples with a notch at $6 \mathrm{mm}$ height (experiments 'Exp 1' to 'Exp 4' versus simulations).}
\label{6mm_force_displ}
\end{figure}

In Figure \ref{6mm_force_displ}, the numerical force-displacement curves have a similar course as the experimentally achieved force-displacement curves. Here, not just the most converging results 'Neo Hooke, $3$ ref' and '$150\%$ fit, $3$ ref' from Figure \ref{compare_hole_no_hole_150fit_NeoHooke} are given in Figure \ref{6mm_force_displ}, but also the results from a second test run, where we geometrically prerefined two times the area from $2 \mathrm{mm}$ below the notch until $1 \mathrm{mm}$ above the circular hole in addition to a globally prerefined mesh with cell diameter $h=0.3$. See Figure \ref{6mm_compare} for the finite element mesh used for the dashed and coloured force-displacement curves given in Figure \ref{6mm_force_displ}. The force-displacement curves from the simulations computed on the geometrically prerefined meshes seem to give better results compared to the experiments than the results based on adaptively refined meshes and a fixed bandwidth $\epsilon$. This observed mesh-sensitivity of phase-field fracture models is widely discussed in the literature, see for example in \cite{ambati2015review,mang2020phase,heister2015primal,mang2020mesh,wick2017modified}. Further, one can see that the length scale $\epsilon$, the $\epsilon$-$h$ relation as well as the incremental step size $\delta t$ can have an impact on the accurate shape of the force-displacement curves also for very simple benchmark tests. Further analysis on this topic is beyond the scope of this work and will be addressed in future. \\

\subsection{Crack paths comparison}\label{comparison}
\begin{figure}[htbp!]
\centering
\includegraphics[height=4.5cm]{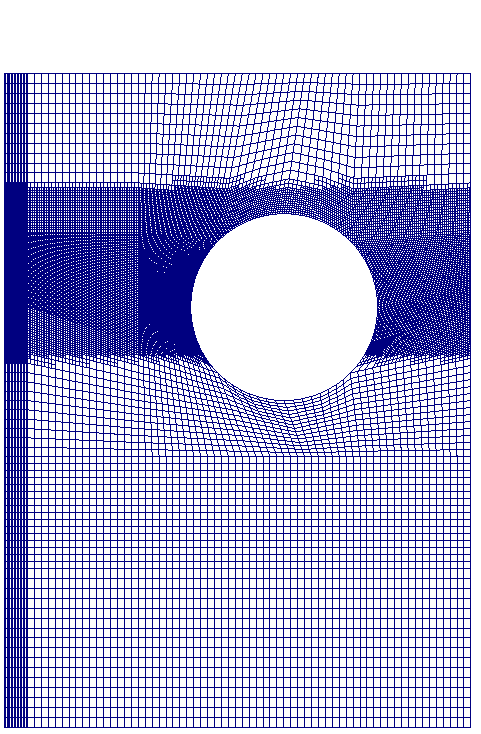}
\includegraphics[height=4.5cm]{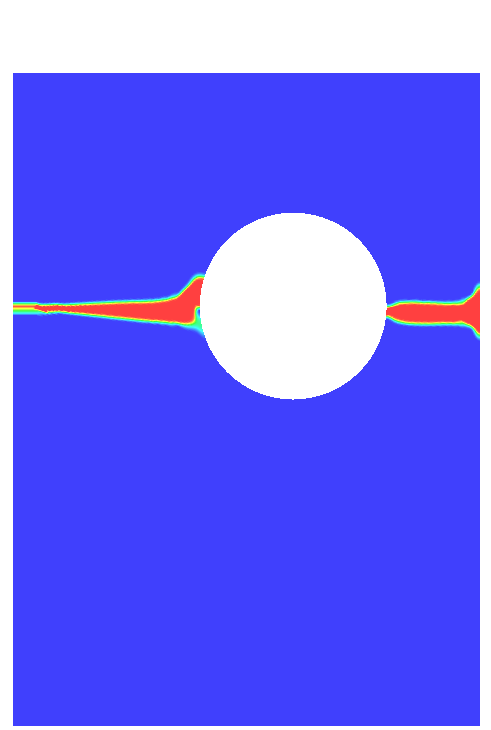}
\includegraphics[height=4.5cm]{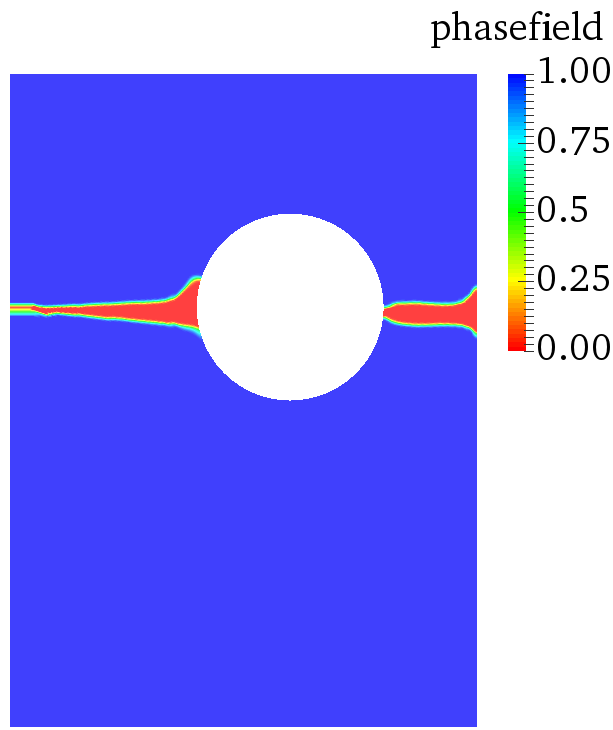}
\hspace{1cm}
\includegraphics[height=4.5cm]{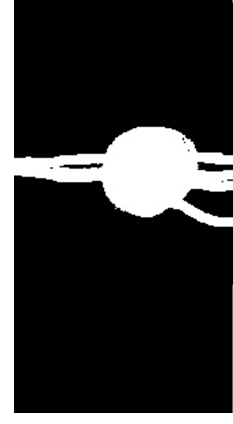}
\caption{Snapshots of the phase-field function after total failure compared to the experimental results for the samples with a notch at $18 \mathrm{mm}$ height. From left to right: the geometrically prerefined mesh, the phase-field function based on the Neo Hookean parameter setting from Table \ref{Tab:Youngsmodulus}, the phase-field function based on the $150\%$ strain fit parameter setting from Table \ref{Tab:Youngsmodulus}, the experimental results from $5$ executed experiments.}
\label{18mm_compare}
\end{figure}

\begin{figure}[htbp!]
\centering
\includegraphics[height=4.5cm]{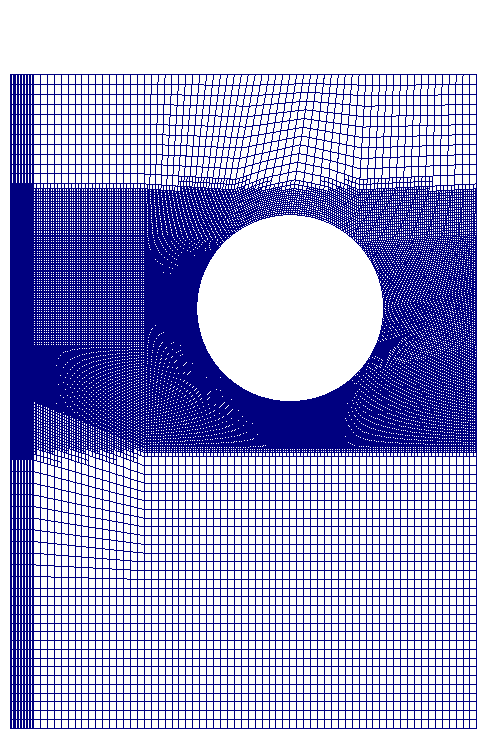}
\includegraphics[height=4.5cm]{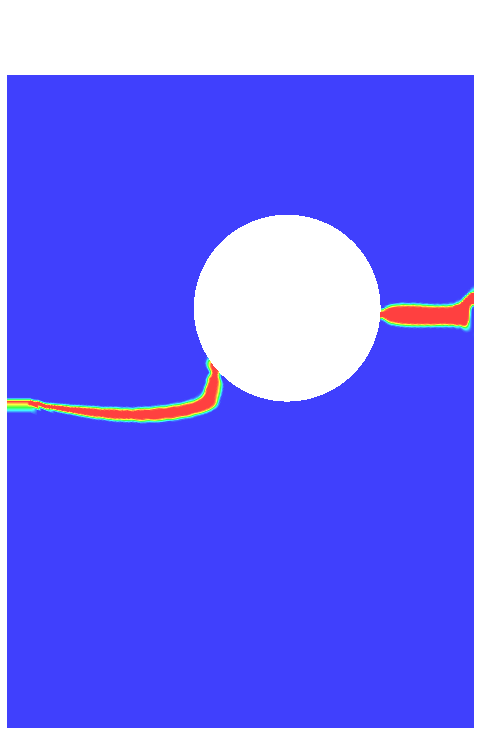}
\includegraphics[height=4.5cm]{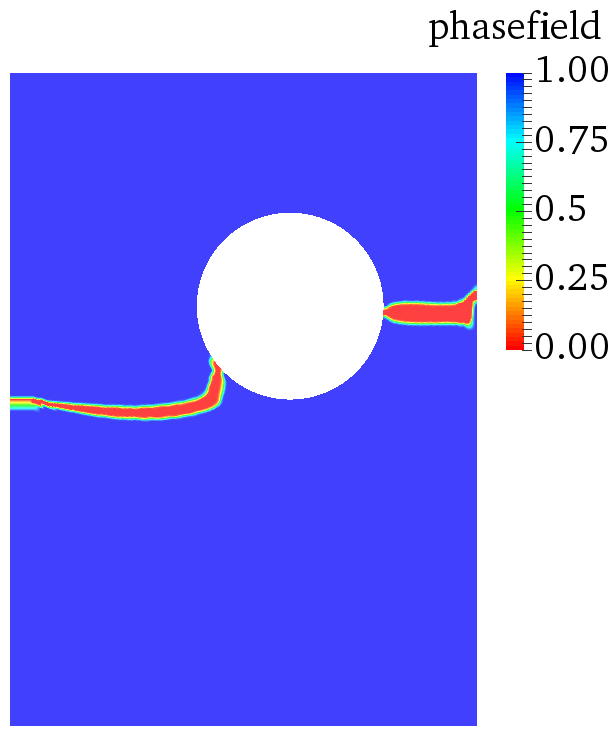}
\hspace{1.15cm}
\includegraphics[height=4.5cm]{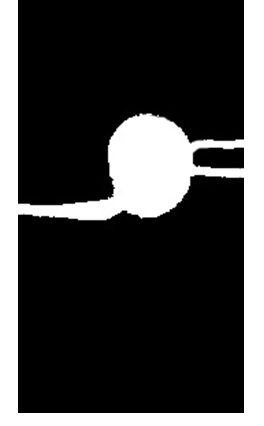}
\caption{Snapshots of the phase-field function after total failure compared to the experimental results for the samples with a notch at $14 \mathrm{mm}$ height. From left to right: the geometrically prerefined mesh, the phase-field function based on the Neo Hookean parameter setting from Table \ref{Tab:Youngsmodulus}, the phase-field function based on the $150\%$ strain fit parameter setting from Table \ref{Tab:Youngsmodulus}, the experimental results from $4$ executed experiments.}
\label{14mm_compare}
\end{figure}

\begin{figure}[htbp!]
\centering
\includegraphics[height=4.5cm]{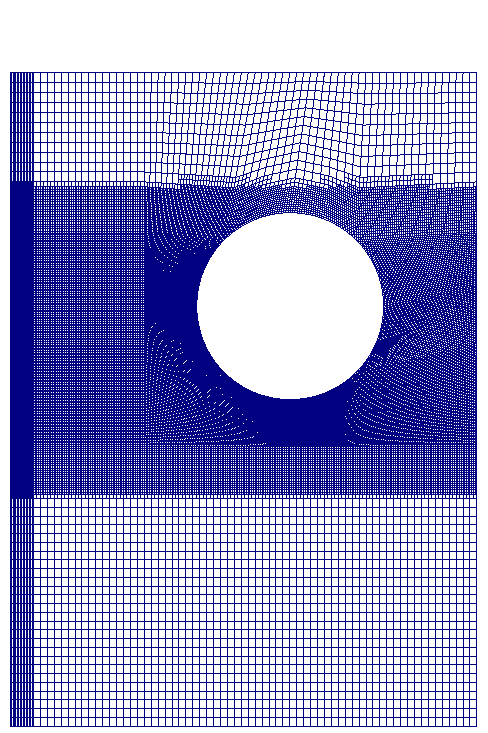}
\includegraphics[height=4.5cm]{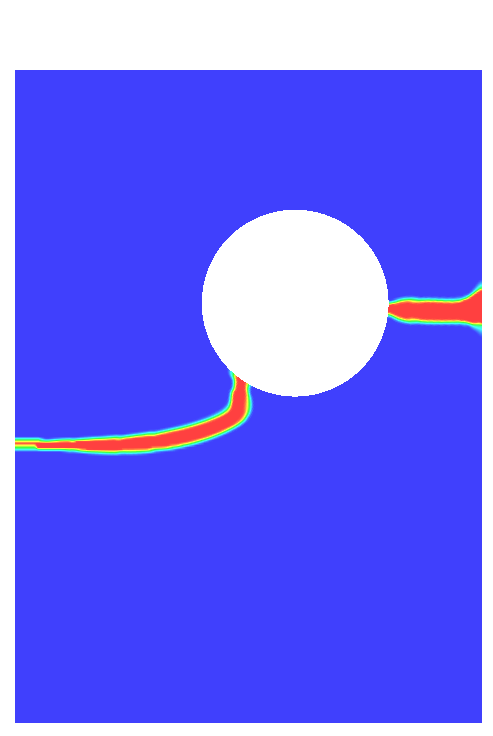}
\includegraphics[height=4.5cm]{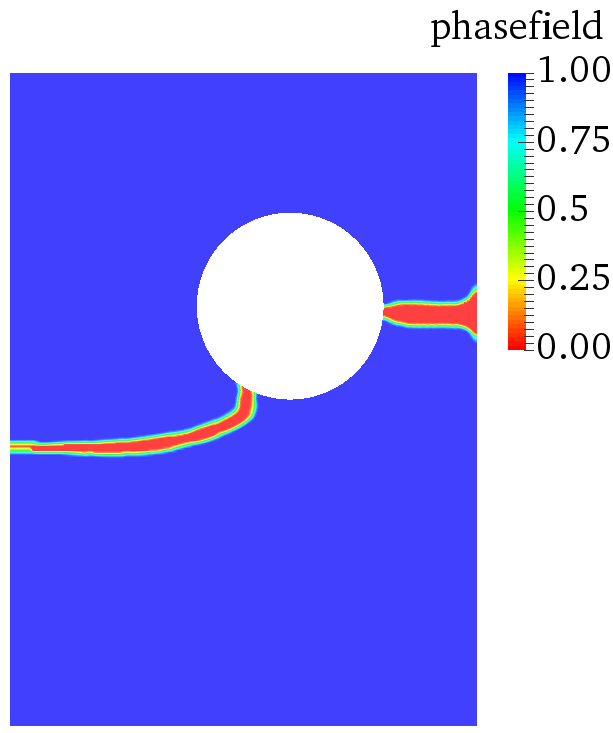}
\hspace{1cm}
\includegraphics[height=4.5cm]{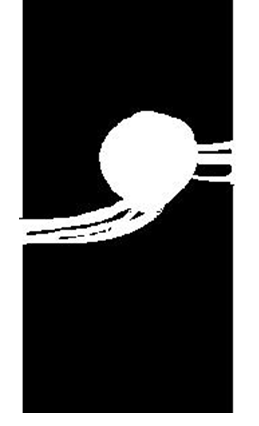}
\caption{Snapshots of the phase-field function after total failure compared to the experimental results for the samples with a notch at $12 \mathrm{mm}$ height. From left to right: the geometrically prerefined mesh, the phase-field function based on the Neo Hookean parameter setting from Table \ref{Tab:Youngsmodulus}, the phase-field function based on the $150\%$ strain fit parameter setting from Table \ref{Tab:Youngsmodulus}, the experimental results from $4$ executed experiments.}
\label{12mm_compare}
\end{figure}

\begin{figure}[htbp!]
\centering
\includegraphics[height=4.5cm]{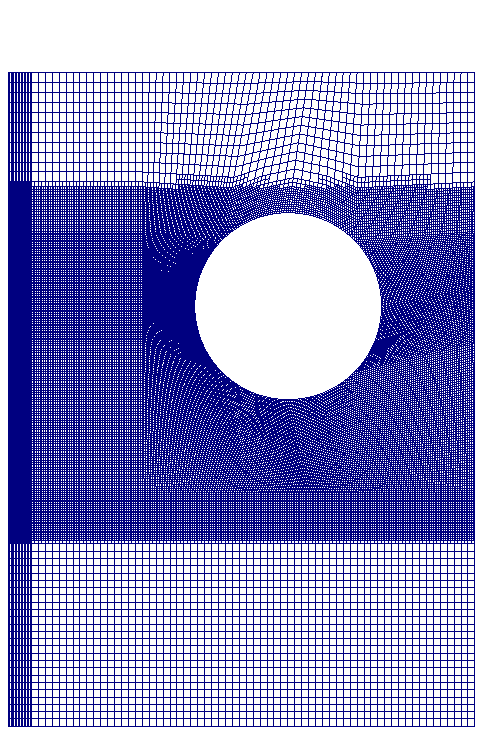}
\includegraphics[height=4.5cm]{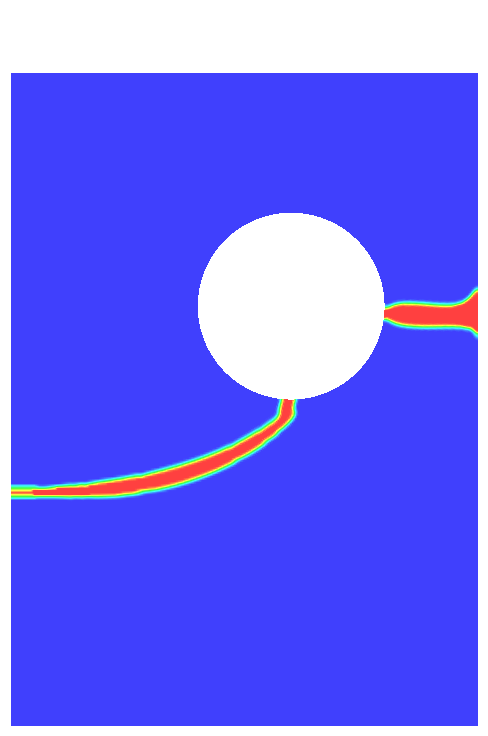}
\includegraphics[height=4.5cm]{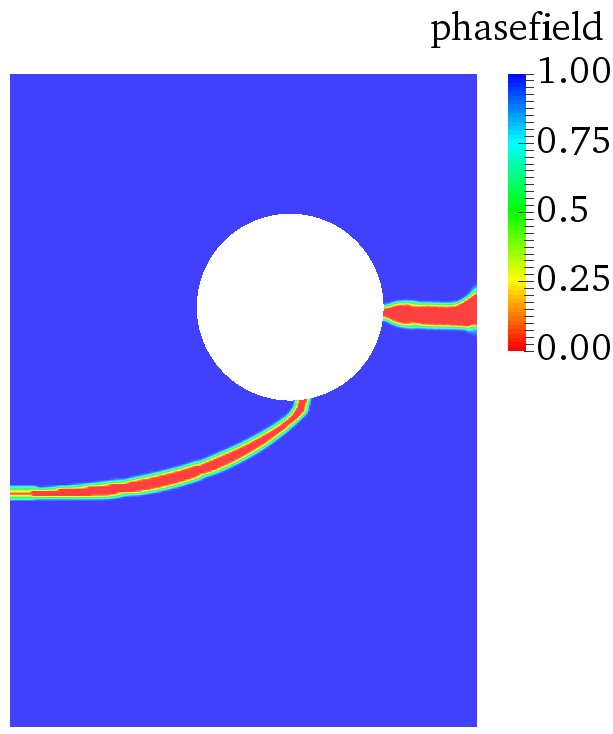}
\hspace{1cm}
\includegraphics[height=4.5cm]{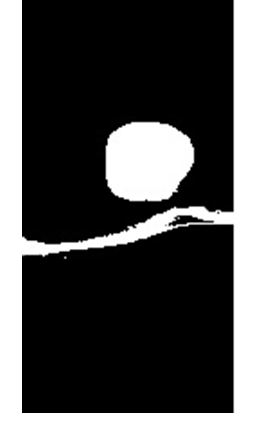}
\caption{Snapshots of the phase-field function after total failure compared to the experimental results for the samples with a notch at $10 \mathrm{mm}$ height. From left to right: the geometrically prerefined mesh, the phase-field function based on the Neo Hookean parameter setting from Table \ref{Tab:Youngsmodulus}, the phase-field function based on the $150\%$ strain fit parameter setting from Table \ref{Tab:Youngsmodulus}, the experimental results from $6$ executed experiments.}
\label{10mm_compare}
\end{figure}

\begin{figure}[htbp!]
\centering
\includegraphics[height=4.5cm]{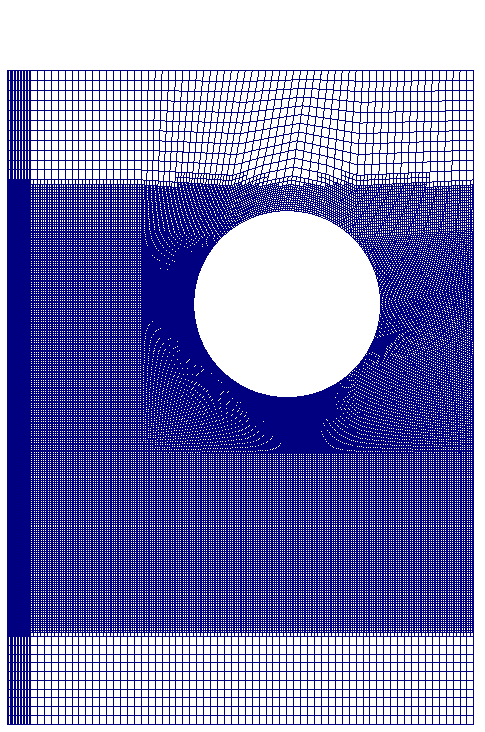}
\includegraphics[height=4.5cm]{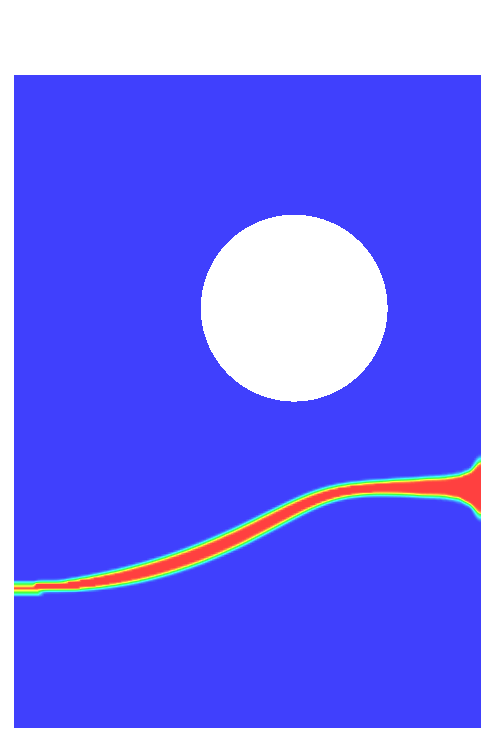}
\includegraphics[height=4.5cm]{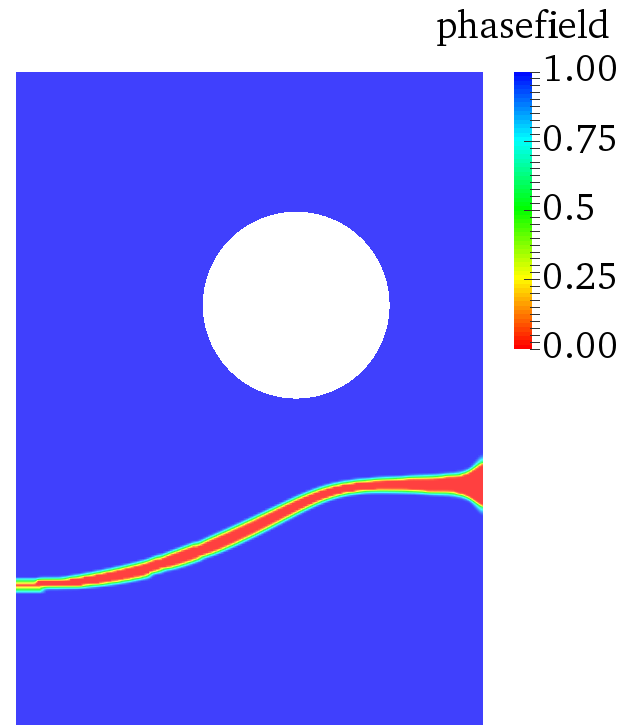}
\hspace{1.3cm}
\includegraphics[height=4.5cm]{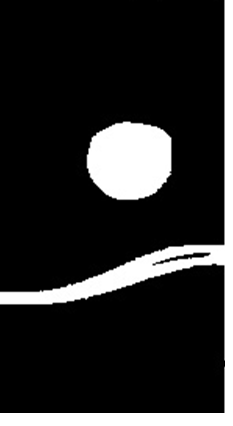}
\caption{Snapshots of the phase-field function after total failure compared to the experimental results for the samples with a notch at $6 \mathrm{mm}$ height. From left to right: the geometrically prerefined mesh, the phase-field function based on the Neo Hookean parameter setting from Table \ref{Tab:Youngsmodulus}, the phase-field based on the $150\%$ strain fit parameter setting from Table \ref{Tab:Youngsmodulus}, the experimental results from $4$ executed experiments.}
\label{6mm_compare}
\end{figure}

In Figure \ref{18mm_compare} to \ref{6mm_compare}, the phase-field functions for five different notch heights are given in comparison to the experimental crack paths after total failure. The finite element meshes are geometrically prerefined generously including the hole and the given notch. The used numbers of degrees of freedom in Figure \ref{18mm_compare} to \ref{6mm_compare} are listed in the lower part of Table \ref{dofs}. 
For the tests with a given notch at $18\mathrm{mm}$, $14\mathrm{mm}$, $12\mathrm{mm}$ and $6\mathrm{mm}$, the crack paths of the simulation coincide sufficiently with the average crack paths from $4$ to $6$ conducted experiments in the Figures \ref{18mm_compare}, \ref{14mm_compare}, \ref{12mm_compare} and \ref{6mm_compare}. Furthermore, the numerically achieved crack paths for the two different parameter settings 'Neo Hooke' and '$150 \%$ fit' are very similar for all five tests with different notch heights. Just the location and the angle from where the crack passes the circular inclusion varies slightly, see for example Figure \ref{12mm_compare}. In Figure \ref{10mm_compare}, snapshots of the phase-field function for both parameter settings are given for the test with a notch at $10\mathrm{mm}$. Observing the experiments for this test case, the crack propagates close to hole but not into the hole. Compare also Figure \ref{fig:NHK00169_v2} for the $10\mathrm{mm}$ test (second row from the bottom), especially the first, third and fourth test specimen, where the crack propagates very close to the hole. The phase-field fracture simulation shows a fracture path similar to them of the $14$ and $18\mathrm{mm}$ tests. One reason could be, that the crack width $\epsilon$ is, from a numerical point of view, not small enough to allow cracks very close to the boundary of the hole without cracking into it.

In the numerical results in Figure \ref{18mm_compare} to \ref{10mm_compare}, we observe that the location, where the crack propagates from the hole to the right boundary of the EPDM strips, is in the very middle of the hole and goes straight to the right (shortest way for the crack) while in the experiments the angles and paths of the crack paths on the right of the hole vary a lot, cf. Figures \ref{18mm_compare}.

\subsection{Comment on our choice of energy functional and strain energy splitting}\label{comment}
Our choice of the energy functional design as well the strain energy splitting in Section \ref{NewModel} is deliberately fitted to the experiments in Section \ref{Experi}.
We made the experience that even with a given notch, the crack starts propagating from the circular hole in the inner of the EPDM strips and propagates simultaneously to the left and right, if the mixed phase-field fracture model is based on an Ambrosio-Tortorelli functional \cite{ambrosio1990approximation,ambrosio1992approximation} (AT$_1$ or AT$_2$, see \cite{tanne2018crack} for a comparison of AT$_1$ and AT$_2$ considering crack nucleation) and the elastic energy splitting approach of Miehe \cite{miehe2010phase}. This high sensitivity around the inclusion also for a very small incremental size $\delta t$ and proper adaptive refinement schemes led us to the energy functional of Wu \cite{wu2017unified} and Amor's volumetric-deviatoric energy splitting \cite{amor2009regularized}. In the literature, for example in \cite{wu2018phase}, a pre-cracked sample with two holes is presented in a compressible solid. Beside the fact, that depending on the mesh, the boundary conditions, initial conditions and the underlying phase-field model, very different crack path results are presented. Further, they point out the need of further studies on problems with holes and inclusions in in \cite[Section 8.2]{wu2018phase}. Other references on computing and discussing crack propagation in solids with inclusions are \cite{zheng2015complementarity,artina2015anisotropic,Loew2019b,kumar2020revisiting}. The last two cited groups are also considering (nearly) incompressible materials.


\section{Discussion of experimental and numerical results}\label{Discussion} 
Despite the assumed simplifications related to geometrical and material non-linearities, the model is able to capture the behavior at least qualitatively.
The experimental data related to the tracking of the crack paths as well as the respective numerical simulations suggest a clear dependency on the position of the initial crack position, which can be confirmed by the numerical simulations. 
In the bottom row in Figure \ref{fig:NHK00169_v2} as well as in Figure \ref{6mm_compare} one can see that even for an initial crack near the lower boundary ($6\mathrm{mm}$) 
far from the hole, the influence of the inhomogeneous stress/strain field is still apparent. The crack starts perpendicularly to the direction of tension on the right side (here: back side), but is diverted upwards later on and returns to a perpendicular path in the end. \\
This effect is observed in the experimental data till initial cracks with a notch height of $10 \mathrm{mm}$. At the end of the crack path for specimens marked with $\#$7 (second row from below), the path is diverted downwards. In the finite element simulation, the crack path within the $10 \mathrm{mm}$ test differs such that the crack propagates into the hole as also seen in the $12, 14$ and $18 \mathrm{mm}$ tests. 
For initial cracks at $12$ and $14 \mathrm{mm}$ height, the path proceeds in a curved shape into the hole and continues perpendicular on the other side of it. 
At $18 \mathrm{mm}$ height, approximately on position of the circle's center, 
the initial crack propagates directly perpendicular into the hole and propagates likewise on the other side until complete rupture. The same can be observed in the numerical results for the tests with a notch height of $12, 14$ and $18 \mathrm{mm}$. Evaluating the statistics of the experiments, repeatability is given. The repetitions of the tests indicate that minor variations, 
e.g.\,small changes in the initial crack length or its angle will have small effects on the general crack path as well. Although dealing with a carbon black filled rubber, no bifurcations in the crack paths for the EPDM compound were observed. Due to the high variation in the onset of total failure in the experiments the calibration of related phase-field parameters might be tricky in the future.


\section{Conclusions}\label{Conclusion}
In this work, a phase-field fracture model in mixed form is proposed to simulate crack propagation in EPDM rubber considering the crack path behavior. The model is based on a mixed form of the solid equation derived from Wu's energy functional \cite{wu2017unified} and the anisotropic model of Amor et al. \cite{amor2009regularized} with a volumetric-deviatoric strain energy splitting, matching the volume conservation of incompressible materials. The obtained crack paths in the numerical results compared to the punctured EPDM strips after stretching them with a certain load until total failure, coincide aside from experimental scattering and the inaccuracy in phase-field fracture modeling with respect to the smeared transition zone around the crack. Even if a quasi-static phase-field fracture model assuming linear elasticity is used for simulating crack propagation in punctured EPDM strips with a given notch, the crack paths and force-displacement results are promising. Especially the crack paths coincide satisfactorily for the tests with a notch height of $18, 14, 12$ and $6 \mathrm{mm}$, see Figures \ref{18mm_compare}, \ref{14mm_compare}, \ref{12mm_compare} and \ref{6mm_compare}, respectively. Furthermore, the force-displacement curves, as an example for the $6\mathrm{mm}$ test, are similar to the experimentally achieved curves, especially with respect to the maximum force at the crack start (Figure \ref{fig:crackmaxforce}) and the traverse displacement at the maximum force (Figure \ref{fig:crackstart}) within the natural scattering of the experiments.
However, we stress that most of our comparisons are of qualitative nature because of simplifications in the governing material models such as nonlinear behavior. Due to the complexity of the experimental and numerical setups, our results are nonetheless a major advancement and themselves novel in the published literature.\\
In future work, a dynamic phase-field fracture model will be considered. Further, a detailed analysis on the mesh sensitivity and $\epsilon$-$h$ relation in simulating crack propagation in rubber-like materials will be addressed. In addition, considering the high non-linear behavior of filled rubbers in cyclic applications, non-linear geometrical effects as well as non-linear material effects have to be addressed. Further, in cyclic dynamic loadings viscoelasticity, temperature effects due to energy dissipation, static hysteresis and permanent deformations become more relevant. Several models to describe those effects are benchmarked in \cite{Carleo2018}, see also \cite{Plagge2020}. Perspectively, the combination of sophisticated mechanical material models and phase field fracture models for describing the crack propagation in rubbers seems promising.

\section*{Acknowledgments}
The authors would like to thank Pénélope Barbery (student at ENSTA Bretagne, Brest) for her support in conducting some of the experiments during her internship at DIK.
Further, the work has been supported by the German Research Foundation, Priority Program 1748 (ID 392587580, WI 4367/2-1).


\bibliography{Literature}

\begin{thebibliography}{10}
\expandafter\ifx\csname url\endcsname\relax
  \def\url#1{\texttt{#1}}\fi
\expandafter\ifx\csname urlprefix\endcsname\relax\def\urlprefix{URL }\fi
\expandafter\ifx\csname href\endcsname\relax
  \def\href#1#2{#2} \def\path#1{#1}\fi

\bibitem{Wohler.1870}
A.~W{\"o}hler, {{\"U}ber die Festigkeitsversuche mit Eisen und Stahl},
  Zeitschrift f{\"u}r Bauwesen 20 (1870) 73--106.

\bibitem{Gehrmann.2019}
O.~Gehrmann, N.~H. Kr{\"o}ger, M.~Krause, D.~Juhre, Dissipated energy density
  as fatigue criterion for non-relaxing tensional loadings of non-crystallizing
  elastomers?, Polymer Testing 78 (2019) 105953 (online).

\bibitem{Ludwig.2017}
M.~Ludwig, {Entwicklung eines Lebensdauer-Vorhersagekonzepts f\"ur
  Elastomerwerkstoffe unter Ber\"ucksichtigung der Fehlstellenstatistik},
  {Dissertation, Universit\"at Hannover (2017)}.

\bibitem{ElYaagoubi2018}
M.~El~Yaagoubi, D.~Juhre, J.~Meier, N.~H. Kr{\"o}ger, T.~Alshuth, U.~Giese,
  Lifetime prediction of filled elastomers based on particle distribution and
  the j-integral evaluation, International Journal of Fatigue 112 (2018)
  341--354.

\bibitem{Gehrmann.2019a}
O.~Gehrmann, M.~El~Yaagoubi, H.~El~Maanaoui, J.~Meier, Lifetime prediction of
  simple shear loaded filled elastomers based on the probability distribution
  of particles, Polymer Testing 75 (2019) 229--236.

\bibitem{Lemaitre.1985}
J.~{Lemaitre}, A continuous damage mechanics model for ductile fracture,
  Journal of Engineering Materials and Technology 107~(1) (1985) 83--89.

\bibitem{Grandcoin.2014}
J.~Grandcoin, A.~Boukamel, S.~Lejeunes, A micro-mechanically based continuum
  damage model for fatigue life prediction of filled rubbers, International
  Journal of Solids and Structures 51~(6) (2014) 1274--1286.

\bibitem{Charrier2003}
P.~Charrier, E.~Ostoja-Kuczynski, E.~Verron, G.~Marckmann, L.~Gornet,
  G.~Chagnon, Theoretical and numerical limitations for the simulation of crack
  propagation in natural rubber components, Constitutive Models for Rubber III
  (2003) 3--10.

\bibitem{Timbrell2003}
C.~Timbrell, M.~Wiehahn, G.~Cook, A.~H. Muhr, Simulation of crack propagation
  in rubber, Constitutive Models for Rubber III (2003) 11--20.

\bibitem{Kaliske2014}
M.~Kaliske, R.~Behnke, R.~Fleischhauer, K.~\"Ozenc, I.~M. Zreid,
  Theretical-numerical approaches to simulate fracture in polymeric materials,
  Procedia Materials Science 3 (2014) 2065--2070.

\bibitem{griffith1920phenomena}
A.~Griffith, The phenomena of flow and rupture in solids, Transactions of the
  Royal Society A 221 (1920) 163--198.

\bibitem{FraMar98}
G.~Francfort, J.-J. Marigo, Revisiting brittle fracture as an energy
  minimization problem, Journal of the Mechanics and Physics of Solids 46~(8)
  (1998) 1319--1342.

\bibitem{bourdin2008variational}
B.~Bourdin, G.~A. Francfort, J.-J. Marigo, The variational approach to
  fracture, Journal of elasticity 91~(1-3) (2008) 5--148.

\bibitem{bourdin2000numerical}
B.~Bourdin, G.~A. Francfort, J.-J. Marigo, Numerical experiments in revisited
  brittle fracture, Journal of the Mechanics and Physics of Solids 48~(4)
  (2000) 797--826.

\bibitem{ambrosio1990approximation}
L.~Ambrosio, V.~M. Tortorelli, Approximation of functional depending on jumps
  by elliptic functional via t-convergence, Communications on Pure and Applied
  Mathematics 43~(8) (1990) 999--1036.

\bibitem{wu2017unified}
J.-Y. Wu, A unified phase-field theory for the mechanics of damage and
  quasi-brittle failure, Journal of the Mechanics and Physics of Solids 103
  (2017) 72--99.

\bibitem{wu2018geometrically}
J.-Y. Wu, A geometrically regularized gradient-damage model with energetic
  equivalence, Computer Methods in Applied Mechanics and Engineering 328 (2018)
  612--637.

\bibitem{wu2018length}
J.-Y. Wu, V.~P. Nguyen, A length scale insensitive phase-field damage model for
  brittle fracture, Journal of the Mechanics and Physics of Solids 119 (2018)
  20--42.

\bibitem{amor2009regularized}
H.~Amor, J.-J. Marigo, C.~Maurini, Regularized formulation of the variational
  brittle fracture with unilateral contact: Numerical experiments, Journal of
  the Mechanics and Physics of Solids 57~(8) (2009) 1209--1229.

\bibitem{ambati2015review}
M.~Ambati, T.~Gerasimov, L.~De~Lorenzis, A review on phase-field models of
  brittle fracture and a new fast hybrid formulation, Computational Mechanics
  55~(2) (2015) 383--405.

\bibitem{wu2018phase}
J.-Y. Wu, V.~P. Nguyen, C.~T. Nguyen, D.~Sutula, S.~Bordas, S.~Sinaie, Phase
  field modeling of fracture, Advances in applied mechancis: multi-scale theory
  and computation 52.

\bibitem{bourdin2019past}
B.~Bourdin, G.~A. Francfort, Past and present of variational fracture, SIAM
  News 52~(9).

\bibitem{wick2020multiphysics}
T.~Wick, Multiphysics Phase-Field Fracture: Modeling, Adaptive Discretizations,
  and Solvers, Vol.~28, Walter de Gruyter GmbH \& Co KG, 2020.

\bibitem{mang2020phase}
K.~Mang, T.~Wick, W.~Wollner, A phase-field model for fractures in nearly
  incompressible solids, Computational Mechanics 65~(1) (2020) 61--78.

\bibitem{mang2020pamm}
A.~Fehse, N.~H. Kr\"{o}ger, K.~Mang, T.~Wick, Crack path comparisons of a mixed
  phase-field fracture model and experiments in punctured epdm strips, PAMM
  20~(1) (2020) pamm.202000335.

\bibitem{Loew2019a}
P.~J. Loew, B.~Peters, L.~A.~A. Beex, Fatigue phase-field damage modeling of
  rubber, Constitutive Models for Rubber XI (2019) 408--412.

\bibitem{Loew2019b}
P.~J. Loew, B.~Peters, L.~A.~A. Beex, Rate-dependent phase-field damage
  modeling of rubber and its experimental parameter identification, Journal of
  the Mechanics and Physics of Solids 127 (2019) 266--294.

\bibitem{Faye2019}
A.~Faye, Y.~Lev, K.~Y. Volokh, Modeling dynamic fracture in rubberlike
  materials, Constitutive Models for Rubber XI (2019) 505--511.

\bibitem{Volohk2017}
K.~Y. Volokh, Fracture as a material sink, Materials Theory 1~(3) (2017) 1--9.

\bibitem{Gehrmann2017}
O.~Gehrmann, N.~H. Kr{\"o}ger, P.~Erren, D.~Juhre, Estimation of the
  compression modulus of a technical rubber via cyclic volumetric compression
  tests, Technische Mechanik 37~(1) (2017) 28--36.

\bibitem{Ricker2019}
A.~Ricker, N.~H. Kr{\"o}ger, Influence of various curing systems and carbon
  black content on the bulk modulus of {EPDM} rubber, Constitutive Models for
  Rubber XI (2019) 200--205.

\bibitem{Roucou1}
D.~Roucou, J.~Diani, M.~Brieu, J.~Witz, A.~Mbiaskop-Ngassa, Experimental
  investigation of elastomer mode i fracture: An attempt to estimate the
  critical strain energy release rate using sent tests, International Journal
  of Fracture 209 (2018) 163--170.

\bibitem{Roucou2}
D.~Roucou, J.~Diani, M.~Brieu, A.~Mbiaskop-Ngassa, Critical strain energy
  release rate for rubbers: single edge notch tension versus pure shear tests,
  International Journal of Fracture 216 (2019) 31--39.

\bibitem{Roucou2019E}
D.~Roucou, J.~Diani, M.~Brieu, A.~Mbiakop-Ngassa, Impact of strain-induced
  softening on the fracture of a carbon-black filled styrene butadiene rubber,
  Constitutive Models for Rubber XI (2019) 528--532.

\bibitem{Ozelo2012}
R.~R.~M. Ozelo, P.~Sollero, A.~L.~A. Costa, An alternative technique to
  evaluate crack propagation path in hyperelastic materials, Tire Science and
  Technology TSTCA 40~(1) (2012) 42--58.

\bibitem{ambrosio1992approximation}
L.~Ambrosio, V.~Tortorelli, On the approximation of free discontinuity
  problems, Bollettino dell'Unione Matematica Italiana 6~(1) (1992) 105--123.

\bibitem{mumford1989optimal}
D.~B. Mumford, J.~Shah, Optimal approximations by piecewise smooth functions
  and associated variational problems, Communications on pure and applied
  mathematics.

\bibitem{tanne2018crack}
E.~Tann{\'e}, T.~Li, B.~Bourdin, J.-J. Marigo, C.~Maurini, Crack nucleation in
  variational phase-field models of brittle fracture, Journal of the Mechanics
  and Physics of Solids 110 (2018) 80--99.

\bibitem{bourdin2014morphogenesis}
B.~Bourdin, J.-J. Marigo, C.~Maurini, P.~Sicsic, Morphogenesis and propagation
  of complex cracks induced by thermal shocks, Physical review letters 112~(1)
  (2014) 014301.

\bibitem{basava2020adaptive}
S.~Basava, K.~Mang, M.~Walloth, T.~Wick, W.~Wollner, Adaptive and
  pressure-robust discretization of incompressible pressure-driven phase-field
  fracture, arXiv preprint arXiv:2006.16566.

\bibitem{mang2020adaptive}
K.~Mang, M.~Walloth, T.~Wick, W.~Wollner, Adaptive numerical simulation of a
  phase-field fracture model in mixed form tested on an l-shaped specimen with
  high poisson ratios, arXiv preprint arXiv:2003.09459.

\bibitem{MieWelHof10a}
C.~Miehe, F.~Welschinger, M.~Hofacker, Thermodynamically consistent phase-field
  models of fracture: variational principles and multi-field fe
  implementations, International Journal for Numerical Methods in Fluids 83
  (2010) 1273--1311.

\bibitem{wick2017error}
T.~Wick, An error-oriented newton/inexact augmented lagrangian approach for
  fully monolithic phase-field fracture propagation, SIAM Journal on Scientific
  Computing 39~(4) (2017) B589--B617.

\bibitem{heister2015primal}
T.~Heister, M.~F. Wheeler, T.~Wick, A primal-dual active set method and
  predictor-corrector mesh adaptivity for computing fracture propagation using
  a phase-field approach, Computer Methods in Applied Mechanics and Engineering
  290 (2015) 466--495.

\bibitem{GeLo16}
T.~Gerasimov, L.~D. Lorenzis, A line search assisted monolithic approach for
  phase-field computing of brittle fracture, Computer Methods in Applied
  Mechanics and Engineering 312 (2016) 276 -- 303.

\bibitem{wick2017modified}
T.~Wick, Modified newton methods for solving fully monolithic phase-field
  quasi-static brittle fracture propagation, Computer Methods in Applied
  Mechanics and Engineering 325 (2017) 577--611.

\bibitem{KoKr20}
A.~Kopanicakova, R.~Krause, A recursive multilevel trust region method with
  application to fully monolithic phase-field models of brittle fracture,
  Computer Methods in Applied Mechanics and Engineering 360 (2020) 112720.

\bibitem{KRISTENSEN2020102446}
P.~K. Kristensen, E.~Martinez-Paneda, Phase field fracture modelling using
  quasi-newton methods and a new adaptive step scheme, Theoretical and Applied
  Fracture Mechanics 107 (2020) 102446.

\bibitem{WAMBACQ2021113612}
J.~Wambacq, J.~Ulloa, G.~Lombaert, S.~François, Interior-point methods for the
  phase-field approach to brittle and ductile fracture, Computer Methods in
  Applied Mechanics and Engineering 375 (2021) 113612.

\bibitem{WU2020102440}
J.-Y. Wu, Y.~Huang, Comprehensive implementations of phase-field damage models
  in abaqus, Theoretical and Applied Fracture Mechanics 106 (2020) 102440.

\bibitem{WU2020112704}
J.-Y. Wu, Y.~Huang, V.~P. Nguyen, On the bfgs monolithic algorithm for the
  unified phase field damage theory, Computer Methods in Applied Mechanics and
  Engineering 360 (2020) 112704.

\bibitem{JoLaWi20}
D.~Jodlbauer, U.~Langer, T.~Wick, Matrix-free multigrid solvers for phase-field
  fracture problems, Computer Methods in Applied Mechanics and Engineering 372
  (2020) 113431.

\bibitem{heister2018parallel}
T.~Heister, T.~Wick, Parallel solution, adaptivity, computational convergence,
  and open-source code of 2d and 3d pressurized phase-field fracture problems,
  PAMM 18~(1) (2018) e201800353.

\bibitem{HeiWi20}
T.~Heister, T.~Wick, pfm-cracks: A parallel-adaptive framework for phase-field
  fracture propagation, Software Impacts 6 (2020) 100045.

\bibitem{arndt2020deal}
D.~Arndt, W.~Bangerth, B.~Blais, T.~C. Clevenger, M.~Fehling, A.~V. Grayver,
  T.~Heister, L.~Heltai, M.~Kronbichler, M.~Maier, et~al., The deal. ii
  library, version 9.2, Journal of Numerical Mathematics 1~(ahead-of-print).

\bibitem{mang2020mesh}
K.~Mang, M.~Walloth, T.~Wick, W.~Wollner, Mesh adaptivity for quasi-static
  phase-field fractures based on a residual-type a posteriori error estimator,
  GAMM-Mitteilungen 43~(1) (2020) e202000003.

\bibitem{miehe2010phase}
C.~Miehe, M.~Hofacker, F.~Welschinger, A phase field model for rate-independent
  crack propagation: {R}obust algorithmic implementation based on operator
  splits, Computer Methods in Applied Mechanics and Engineering 199 (2010)
  2765--2778.

\bibitem{zheng2015complementarity}
H.~Zheng, F.~Liu, X.~Du, Complementarity problem arising from static growth of
  multiple cracks and mls-based numerical manifold method, Computer Methods in
  Applied Mechanics and Engineering 295 (2015) 150--171.

\bibitem{artina2015anisotropic}
M.~Artina, M.~Fornasier, S.~Micheletti, S.~Perotto, Anisotropic mesh adaptation
  for crack detection in brittle materials, SIAM Journal on Scientific
  Computing 37~(4) (2015) B633--B659.

\bibitem{kumar2020revisiting}
A.~Kumar, B.~Bourdin, G.~A. Francfort, O.~Lopez-Pamies, Revisiting nucleation
  in the phase-field approach to brittle fracture, Journal of the Mechanics and
  Physics of Solids (2020) 104027.

\bibitem{Carleo2018}
F.~Carleo, E.~Barbieri, R.~Whear, J.~Busfield, Limitations of viscoelastic
  constitutive models for carbon-black reinforced rubber in medium dynamic
  strains and medium strain rates, Polymers 10~(9) (2018) 988.

\bibitem{Plagge2020}
J.~Plagge, A.~Ricker, N.~Kröger, P.~Wriggers, M.~Klüppel, Efficient modeling
  of filled rubber assuming stress-induced microscopic restructurization,
  International Journal of Engineering Science 151 (2020) 103291.

\end{thebibliography}

\end{document}